\documentclass[12pt]{article}
\usepackage{amsfonts}
\usepackage{epsfig}
\title{Outer Billiards and The Pinwheel Map}
\author{Richard Evan Schwartz \thanks{\hskip 5 pt Supported by 
N.S.F. Research Grant DMS-0072607}}

\newtheorem{theorem}{Theorem}[section]

\newtheorem{lemma}[theorem]{Lemma}

\newtheorem{corollary}[theorem]{Corollary}

\def\startproof{{\bf {\medskip}{\noindent}Proof: }}

\def\endproof{$\spadesuit$  \newline}

\def\R{\mbox{\boldmath{$R$}}}%
\def\Z{\mbox{\boldmath{$Z$}}}%

\begin{document}
\maketitle

\section{Introduction}

\subsection{Background}

B. H. Neumann [{\bf N\/}] introduced outer billiards in
the late 1950s and J. Moser [{\bf M1\/}]
popularized the system in the 1970s as a toy model
for celestial mechanics.
Outer billiards is a discrete self-map of $\R^2-P$,
where $P$ is a bounded convex planar set as
in Figure 1.1 below.
Given $x_1 \in \R^2-P$, one
defines $x_2$ so that
the segment $\overline{x_1x_2}$ is tangent to $P$ at its
midpoint and $P$ lies to the right of the ray
$\overrightarrow{x_1x_2}$.  The map $x_1 \to x_2$
is called {\it the outer billiards map\/}.
The map is almost everywhere defined and
invertible.  

\begin{center}
\resizebox{!}{2.2in}{\includegraphics{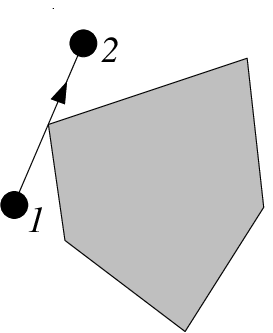}}
\newline
{\bf Figure 1.1:\/} Outer billiards relative to $P$.
\end{center}

The purpose of this paper is to establish
an equivalence between polygonal outer billiards
and an auxilliary map which we call the {\it pinwheel map\/}.
We call our main result the Pinwheel Theorem.
We worked out this equivalence for a restricted version
of outer billiards on kites\footnote{A kite is a convex
quadrilateral having a diagonal that is a line of symmetry.} in
[{\bf S2\/}, Pinwheel Lemma].  

A straightforward version of the Pinwheel Theorem, which
works for all points sufficiently far from the polygon,
appears in almost every paper on polygonal outer billiards.
(See \S \ref{lpm} of this paper for a description and proof.)
Here we mention three papers specifically.
 In Vivaldi-Shaidenko [{\bf VS\/}], Kolodziej [{\bf Ko\/}], and
Gutkin-Simanyi [{\bf GS\/}], it is proved (each with different methods)
that outer billiards on a {\it quasirational polygon\/} has
all orbits bounded.  (See \S \ref{qrdef} for a definition.)

For other work on outer billiards, see
[{\bf M2\/}], [{\bf D\/}] (bounded orbits for
sufficiently smooth convex domains),
[{\bf G\/}] (bounded orbits for trapezoids),
[{\bf S1\/}] (unbounded orbits for the Penrose kite),
[{\bf T2\/}] (aperiodic orbits for the regular pentagon),
[{\bf DF\/}] (unbounded orbits for the half-disk).

In contrast to the ``far away Pinwheel Theorem'', which has
an almost instantaneous proof, the full result is much
more subtle.  The full result
makes a statement about all outer billiards orbits, and
not just those sufficiently far from the polygon.
This stronger statement allows us to give a kind of
bijection between the unbounded orbits of the
pinwheel map and the unbounded outer billiards
orbit -- a bijection that is not deducible just from the
``far away'' result.    See \S \ref{qr0} for precise
statements.
In [{\bf S2\/}] we use this bijection to
show that outer billiards has unbounded
orbits relative to any irrational kite.

Another motivation for studying the pinwheel map
is that it has a geometrically appealing acceleration -- i.e.,
speed up of the time parameter.  This acceleration,
in turn, has a higher dimensional compactification as a polytope
exchange map.   The one case we worked out,
 [{\bf S2\/}, Master Picture Theorem],  had a very rich
structure and was quite decisive for
our theory of outer billiards on kites.  The general case
promises to be equally rich, and we hope that it will be
useful in studying the fundamental questions about
general polygonal billiards -- e.g., the existence of
unbounded orbits.

\subsection{The Main Result}
\label{main result}
\label{qrdef}

We will prove our result for 
convex polygons that have no parallel sides.
We call such polygons {\it nice\/}.   There is
probably a similar result that works for
any convex polygon, but allowing parallel sides
introduces annoying complications.

\begin{center}
\resizebox{!}{2.6in}{\includegraphics{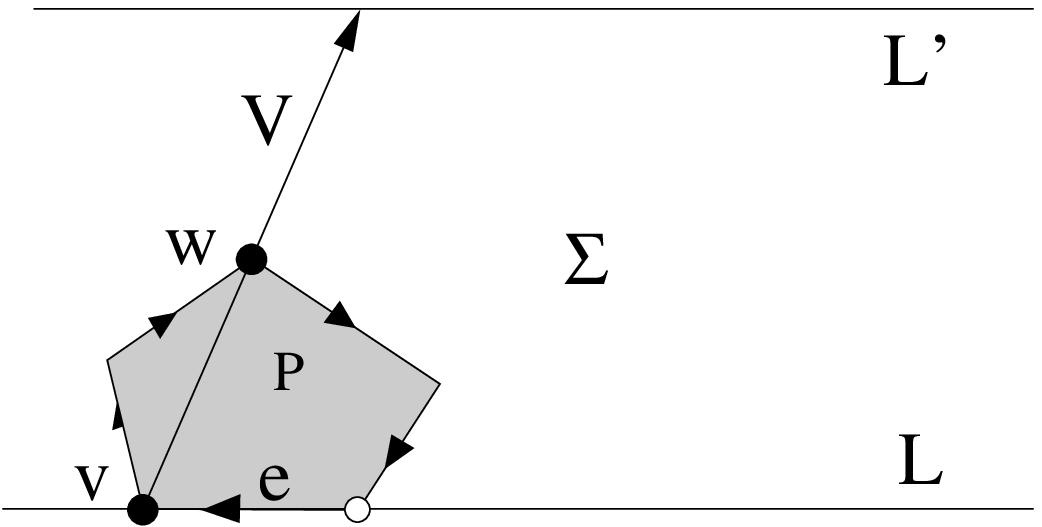}}
\newline
{\bf Figure 1.2:\/} The strip associated to $e$.
\end{center}  

Let $P$ be a nice $n$-gon.  
We orient the edges of $P$ so that
they go clockwise.
 To each edge $e$ of $P$, one
associates a pair $(\Sigma,V)$, where $\Sigma$ is an
infinite strip in the plane and $V$ is a vector that points
from one edge of $\Sigma$ to the other.
Let $v$ be the head vertex of $e$.  Since
$P$ is a nice polygon, there is a unique vertex
$w$ of $P$ that is as far as possible from the
line $L$ containing $e$.  Let $V=2(w-v)$.
 Let $L'$ be the line
parallel to $L$ such that $w$ is equidistant
from $L$ and $L'$. We let $\Sigma$ be the
strip bounded by $L$ and $L'$.  This is  a
construction that comes up often in 
polygonal outer billiards.  We
call $(\Sigma,V)$ a {\it pinwheel pair\/}
and we call $\Sigma$ a {\it pinwheel strip\/}.

The $n$ pinwheel strips $\Sigma_1,...,\Sigma_n$
are cyclically ordered, according to their slopes.
Our convention is that the strips rotate
counterclockwise as we move forward through the
indices.   
\newline
\newline
{\bf Remark:\/}
The $n$-gon $P$ is {\it quasirational\/} if and only if
it may be scaled so that the $n$ parallelograms
$\Sigma_j \cap \Sigma_{j+1}$ all have integer
areas.  Here indices are taken cyclically.
\newline

Given the pair $(\Sigma,V)$, we define a map
$\mu$ on $\R^2-\partial \Sigma$ as follows.
\begin{itemize}
\item If $x \in \Sigma-\partial \Sigma$ then $\mu(x)=x$.
\item If $x \not \in \Sigma$ then $\mu(x)=x \pm V$,
whichever  point is closer to $\Sigma$.
\end{itemize}
The map $\mu$ moves points ``one step'' closer to
lying in $\Sigma$, if they don't already lie in $\Sigma$.
Note that  $\mu$ is not defined on the boundary $\partial \Sigma$.

Let $\R^2_n=\R^2 \times \{1,...,n\}$, with indices taken
mod $n$.  We define the
{\it pinwheel map\/} $\psi^*: \R^2_n \to \R^2_n$ by the
following conditions.
\begin{itemize}
\item $\psi^*(p,k)=(\mu_{k+1}(p),k+1)$ if $\mu_{k+1}(p)=p$.
\item $\psi^*(p,k)=(\mu_{k+1}(p),k)$ if $\mu_{k+1}(p) \not =p$.
\end{itemize}
In other words, we try to move $p$ by the $(k+1)$st strip map.
If the point doesn't move, we increment the index and
give the next strip map a chance to move the point.

Let $\psi$ denote the second iterate of the
outer billiards map. We define $\psi$ to be
the identity inside the polygon $P$. 
Let $\pi: \R_n^2 \to \R^2$ be the projection map.
A {\it section\/} is a map $\iota: \R^2 \to \R^2_n$
such that $\pi \circ \iota$ is the identity.

\begin{theorem}[Pinwheel]
There is a section $\iota: \R^2 \to \R^2_n$ such that
$$\psi(p)=\pi \circ (\psi^*)^k \circ \iota; \hskip 30 pt
k=k(p) \in \{1,...,3n\}.$$
This relation holds on all points for which $\psi$ is
well defined. 
\end{theorem}

Far from the origin, we
have $k(p)=1$ unless $\psi(p)$ lies in a pinwheel
strip.  In this case $k(p)=2$ because one extra
iterate of $\psi^*$ is required to shift the index.
As we mentioned above, the
Pinwheel Theorem is well known, and very easy to
prove for points far from the polygon. (See
Lemma \ref{far}.)  

The new information given by the Pinwheel Theorem is
that correspondence extends in some way to the whole plane.  When
$k(p)>2$ it means that there is a funny cancellation
that happens in order to make the two systems line up.
The fact that $k(p) \leq 3n$ puts a bound on the
complexity of this cancellation.  The bound we
get on $k(p)$ probably isn't sharp but it has
the right order of dependence on $n$.

\subsection{Corollaries}
\label{qr0}

Just from knowledge of the relation far from
the origin, one can conclude nothing about how the
unbounded orbits of one system compare to the unbounded
orbits of the other.  It is easy to construct two
maps of the plane that agree outside a compact set,
such that one of the maps has unbounded orbits and
the other one doesn't.  

The Pinwheel Theorem adds the information that
rules out such wierd pathologies.
We say that an unbounded orbit
of $\psi^*$ is {\it natural\/} if it lies in
$\iota(\R^2)$ sufficiently far from the origin.
In the next result, the word {\it orbit\/} means
both the forwards and backwards orbit. This
makes sense because both $\psi$ and $\psi^*$ are
invertible.  We state our result for the
{\it forward direction\/} of the orbit, though
a similar statement holds for the backwards
direction, and for both directions at the same time.

\begin{corollary}
\label{bijection}
Relative to any nice polygon,
there is a canonical  bijection between the forward unbounded
$\psi$-orbits and the natural forward unbounded
$\psi^*$-orbits.  The bijection sending the orbit
$O$ to the orbit $O^*$ is such that
$O=\pi(O^*)$ outside a compact subset.
\end{corollary}

One reason why one might want to study
$\psi^*$ in place of $\psi$ is that
$\psi^*$ has an appealing acceleration.
Define
\begin{equation}
\label{X}
\widehat X=\bigcup_{j=1}^n (\Sigma_j \times \{j\}) \subset \R_n^2.
\end{equation}
Topologically, $\widehat X$ is the disjoint union of $n$ strips.
$\widehat X$ agrees with $\iota(\R^2)$ outside a compact set.
We let
$\widehat \psi: \widehat X \to \widehat X$ be the first return map of
$\psi^*$ to $X$.  Geometrically, we start with a point
in $\Sigma_1$ and iterate $\mu_2$ until we
land in $\Sigma_2$, then iterate $\mu_3$ until
we land in $\Sigma_3$, etc.  Once again, we state
the result in terms of the forward direction just for
convenience.

\begin{corollary}
\label{bijection2}
There is a canonical bijection between the forward unbounded
orbits of $\psi$ and the forward unbounded orbits of
$\widehat \psi$.
\end{corollary}

One can accelerate somewhat further.
The map $\widehat \Psi = (\widehat \psi)^n$
preserves each individual strip in $\widehat X$.  The
action on each strip is conjugate to the
action on any of the other ones outside of
a compact set.  Thus, we can pick on of the
strips, say $\Sigma_1$, and consider the
map $\widehat \Psi: \Sigma_1 \to \Sigma_1$.
We call $\widehat \Psi$ the {\it pinwheel
return map\/}.

At the same time, we can consider the first return map of
$\psi$ to $\Sigma_1$.  We call this map $\Psi$.  
Again, it is well known that
$\Psi=\widehat \Psi$ outside a compact set, 
but this is not enough information to 
produce a correspondence between the unbounded
orbits of the two systems.
Our next corollary fills in this information.
Once again, we have picked out the forward
direction just for convenience.

\begin{corollary}
\label{bijection3}
There is a canonical bijection between the
forward unbounded orbits
$\Psi$ and the forward unbounded orbits of $\widehat \Psi$.
\end{corollary}

This last corollary is pretty close to a direct
generalization of the Pinwheel Lemma we
proved for kites in [{\bf S2\/}].  In 
[{\bf S2\/}] we actually proved that
$\Psi=\widehat \Psi$ if $\Sigma_1$ was
properly chosen from amongst the
$4$ possible strips.  However, this
stronger result does not seem to be
true in general.

\subsection{Outline of the Paper}

In \S 2 we reduce the Pinwheel Theorem
to two auxilliary results,
Theorem \ref{structure3} and Theorem
\ref{structure2}.  Roughly, Theorem \ref{structure2}
describes an infinite family of ``local substitution
rules'' that allows one to convert between
between the pinwheel map and the polygonal
outer billiards map.  Theorem \ref{structure3}
guarantees that all these local rules patch
together globally.

In \S 3 and \S 4 we
establish some geometric and combinatorial
facts that we will need for the proofs of
Theorems \ref{structure3} and \ref{structure2}.
In \S 5 we prove Theorem \ref{structure3}, and
in \S 6-7 we prove Theorem \ref{structure2}.
Ultimately, the argument boils down to robust
general properties like convexity and induction,
but only after we find the combinatorial
structure of what is going on.  

In \S \ref{qr} we prove
the theorem, due to [{\bf VS\/}], [{\bf K\/}], and
[{\bf GS\/}], that all outer billiards orbits are
bounded relative to a quasi-rational polygon.
Our proof only uses the trivial part of
the Pinwheel Theorem, namely the part
that works for points sufficiently far from the polygon.
There is nothing new in our proof, but we think
that we have boiled the matter down
to its essence.  We include this proof because we already
have the notation set up, and because we like the result so much.

I would like to thank John Smillie and Sergei
Tabachnikov for helpful conversations
about outer billiards and the pinwheel map.
This work is part of an ongoing conversation
John Smillie about the pinwheel map.  In
particular, the formulation we give of
the Pinwheel Theorem emerged in conversations
with Smillie.

I would also like to thank the Clay Mathematics
Institute for their support during my sabbatical
this year.

\newpage

\section{Proof in Broad Strokes}

\subsection{The Forward Partition}

Let $P$ be a nice $n$-gon. 
Let $\psi$ be the square outer billiards map.
We have already mentioned that $\psi$ is a piecewise translation.
For almost every point $p \in \R^2-P$, there is
a pair of vertices $(v_p, w_p)$ of $P$ such that
\begin{equation}
\psi(p)=p+V_p; \hskip 30 pt V_p = 2(w_p-v_p).
\end{equation}
See Figure 2.1.

\begin{center}
\resizebox{!}{2in}{\includegraphics{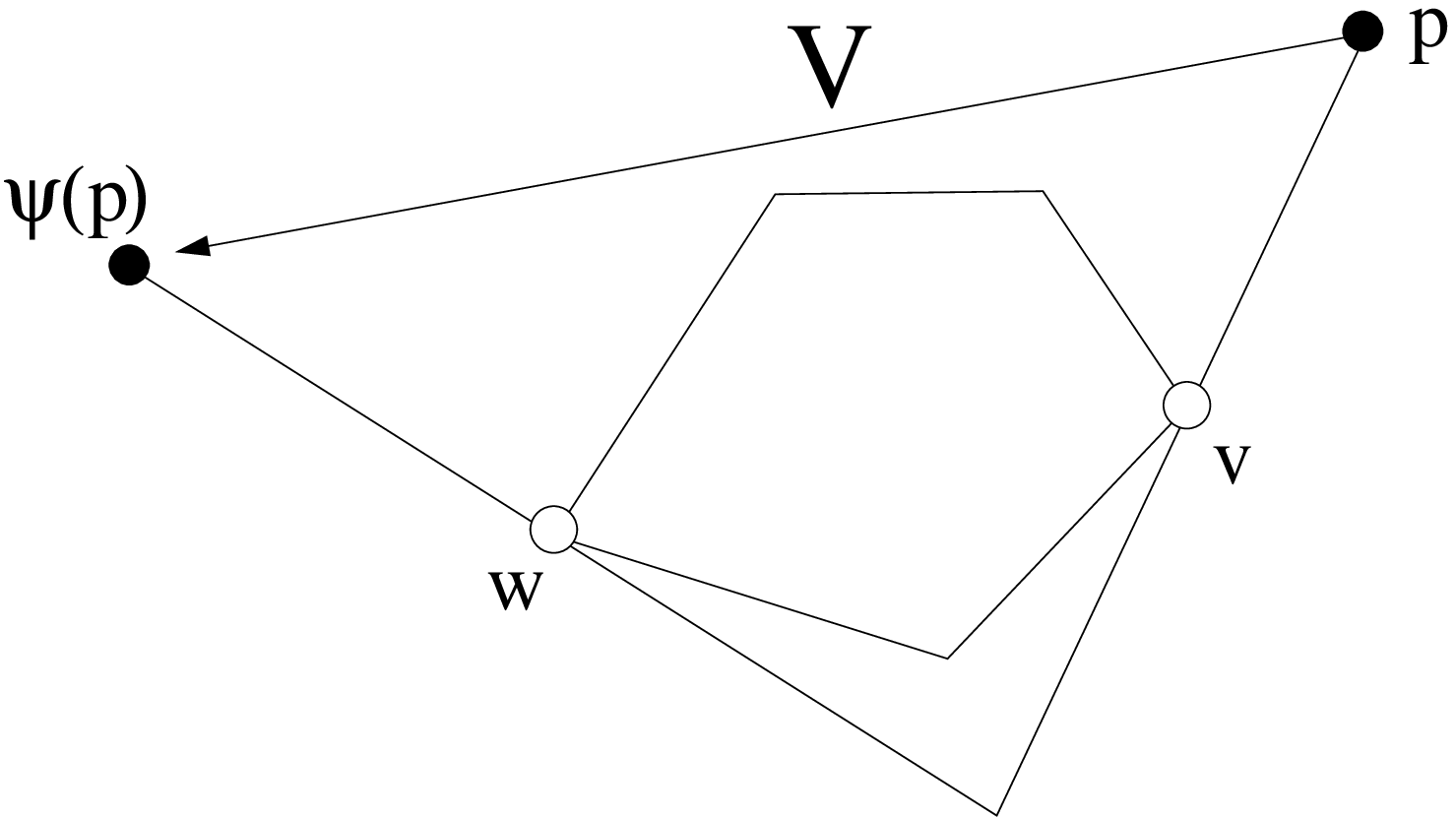}}
\newline
{\bf Figure 2.1:\/} A piecewise translation
\end{center}

The dependence of $V_p$ on $p$ is locally constant,
and the regions where the map $p \to V_p$ is constant
are convex and polygonal.  Since there are only finitely many
pairs of vertices of $P$, we have a partition of
$\R^2-P$ into convex polygonal sets.  We call this
partition the {\it forward partition\/} associated
to $P$.  Here we summarize some of the results we
establish in the next chapter.
\begin{itemize}
\item $P$ has $2n$ unbounded regions, which
we provisionally call $R_i(\pm)$ for $i=1,...,n$.
The region $R_i(\pm)$ is labelled by the pair of vertices
$(v,w)$ such that $2(w-v) = \pm V_i$.  Compare
Figure 1.2.
\item Sufficiently far from the origin, the regions
$R_i(+)$ and $R_i(-)$ share a ray with $\Sigma_i$.
A sufficiently large circle centered at the origin,
oriented counterclockwise, encounters 
$R_1$, $\Sigma_1$, $R_2$, $\Sigma_2$, etc.
The ordering of the plusses and minuses depends
on the geometry of $P$. See Figure 2.3 below.
\end{itemize}

Figure 2.2 shows an example of the forward partition for
a nice octagon.  The central shaded
figure is the octagon. The white regions surrounding
the octagon are the compact tiles of the partition.
The shaded regions on the outside, which have been
cut off by the bounding box of the figure, are
 unbounded tiles.   One of the compact tiles has
also been cut off by the bounding box.

\begin{center}
\resizebox{!}{5.2in}{\includegraphics{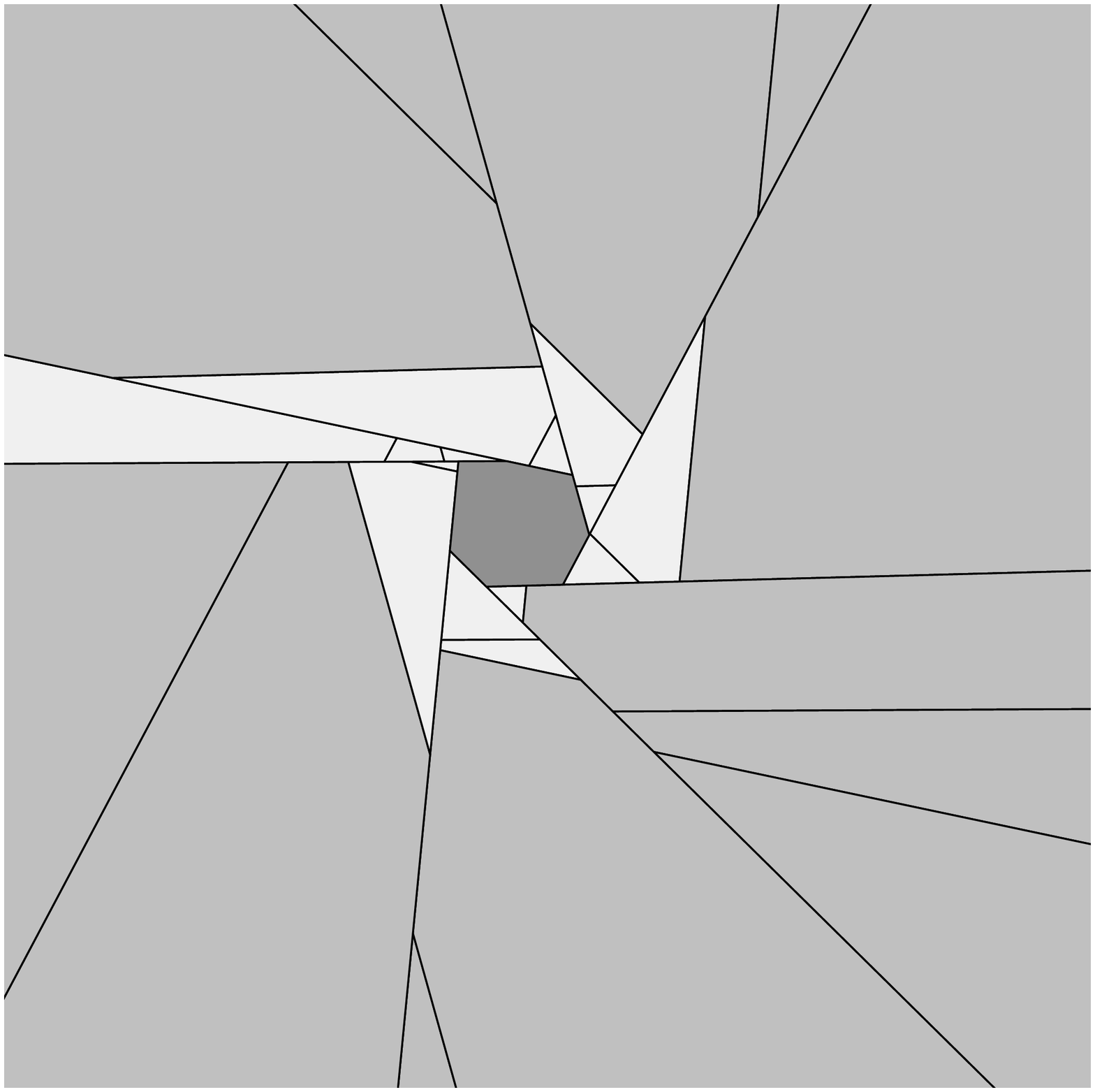}}
\newline
{\bf Figure 2.2:\/} The forward partition for an octagon.
\end{center}

\subsection{Far From the Origin}
\label{lpm}

Let $\mu_1,...,\mu_n$ be the strip maps defined
in \S \ref{main result}.  Let $\widehat X$ be
as in Equation \ref{X}. Given a compact set $K$,
let 
\begin{equation}
\widehat X_K=\widehat X-\pi^{-1}(K); 
\hskip 30 pt X_k=(\Sigma_1 \cup ... \cup \Sigma_n)-K.
\end{equation}
If $K$ is sufficiently large, the map
$\pi: \widehat X_K \to X_K$ is a bijection (and
a local isometry.)
The point is that the pinwheel strips are all
disjoint far from the origin, thanks to the
fact that $P$ is a nice polygon.

\begin{lemma}
\label{far}
$\psi=\pi \circ \psi^* \circ \pi^{-1}$
on $X_K$ provided that $K$ is sufficiently large.
\end{lemma}

\startproof
Figure 2.3 shows the situation for $n=4$.
The central disk $\Delta$ covers up all the
bounded tiles.  The dark lines are parts of
the pinwheel strip boundaries.  We just show the
top half of the picture.

\begin{center}
\resizebox{!}{2.8in}{\includegraphics{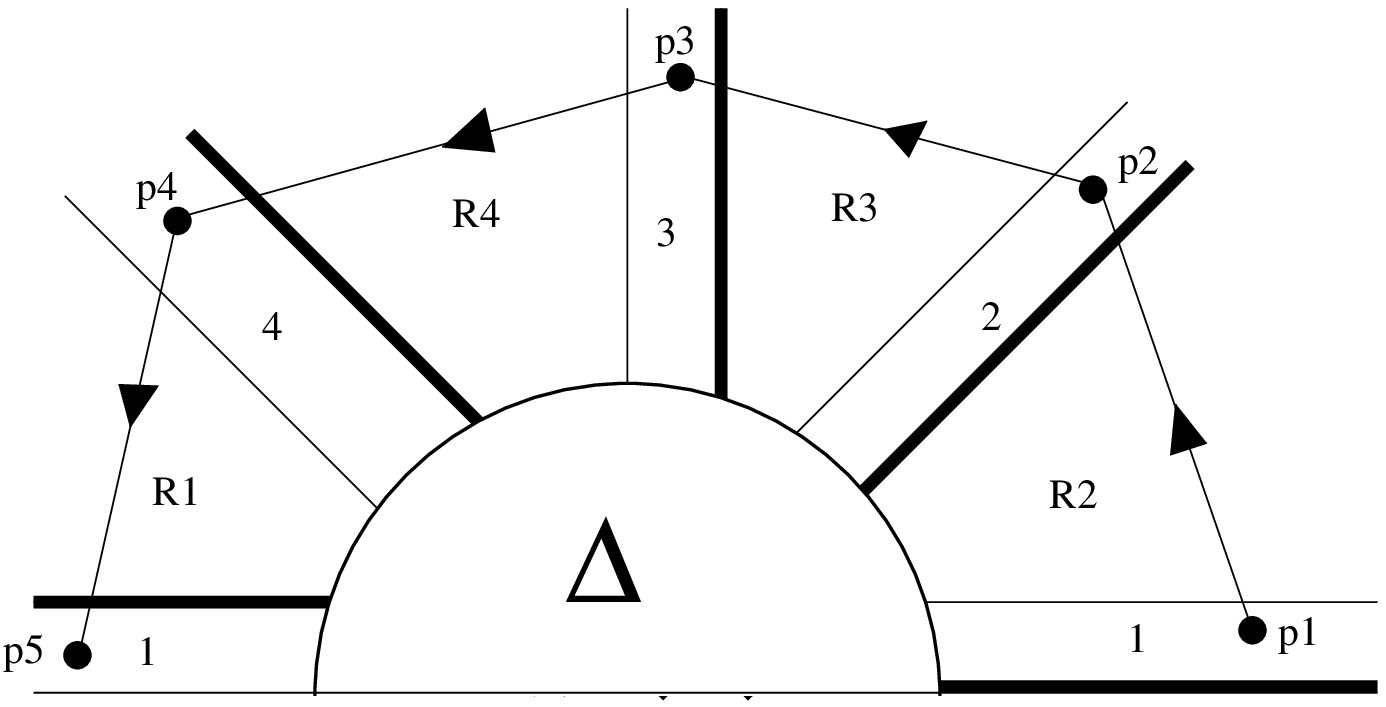}}
\newline
{\bf Figure 2.3:\/} Outside a compact set
\end{center}

\startproof
We will take $k=1$ and set $p_1=p$.
Provided that
$p_1 \in \Sigma_1$ starts sufficiently far away from
$\Delta$, the map $\psi$ simply adds $V_1$ to 
$p_1$ until the resulting orbit reaches a
point $p_2 \in \Sigma_2$.  The map
$\psi^*$ does exactly the same thing, and after
the same number steps, we arrive at the
point $(p_2,2)$.   And so on.
\endproof

\subsection{The Pinwheel Identities}
\label{spoke def}

We have already associated $n$ vectors to
our nice $n$-gon $P$, namely
$V_1,...,V_n$.  Each $V_k$ defines
a line segment $S_k$ that joins the
two ends of $S_k$.  Referring to Figure 1.2,
the spoke $S_k$ joins the vertices $v$ and $w$.
There are $n$ spokes associated to $P$ and
we write these as $S_1,..,S_n$. This cyclic
ordering is compatible with the ordering that
comes from the slopes of the spokes.
Figure 2.4 shows an example, with some of the
edges highlighted.

\begin{center}
\resizebox{!}{3.8in}{\includegraphics{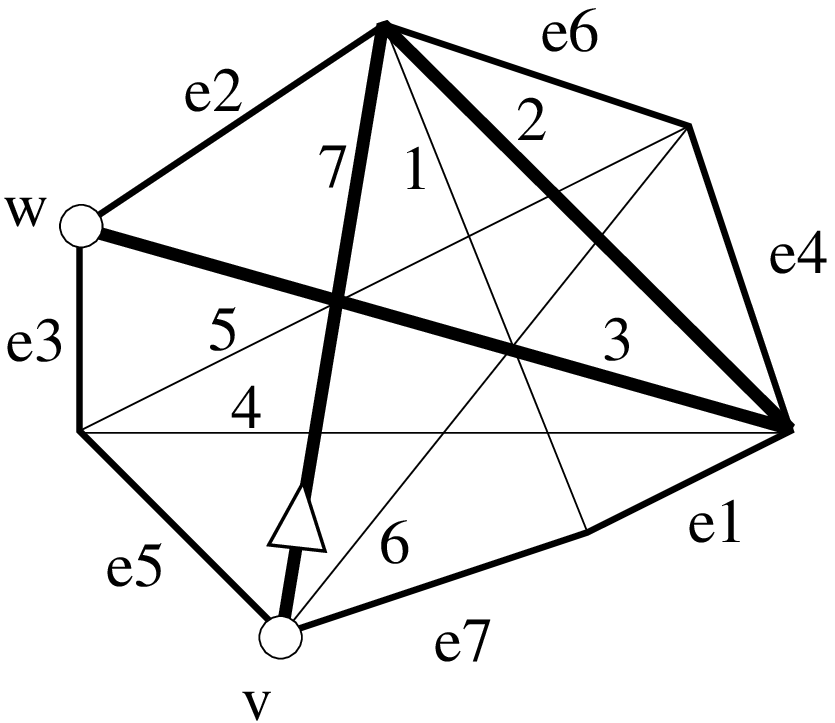}}
\newline
{\bf Figure 2.4:\/} The spokes of the polygon, and the path $7 \to 3$.
\end{center}

We say that an oriented, connected
 polygonal path $\gamma$ is {\it admissible\/} if
the following holds.
\begin{itemize}
\item $\gamma$ consists of an odd number of spokes of $P$.
\item The ordering on the spokes of $\gamma$ is
compatible with the cyclic order.
\item Let $\gamma'$ be the polygonal path in
$\R \cup \infty$ obtained by connecting the slopes
of the spokes of $\gamma$.  Then $\gamma'$ is
a proper subset of $\R \cup \infty$.
\end{itemize}

The last condition means that $\gamma$ does not
wrap all the way around $P$.
We use the notation $a \to b$ to name
admissible paths.  The first spoke is $a$
and the last one is $b$.
We often take
the indices mod $n$ and use $b+n$ in
place of $b$ in case $b<a$.  Thus,
$7 \to 10$ is another name for the
the path in Figure 2.5.

We prove the following result in \S \ref{admissible proof}

\begin{lemma}
\label{structure1}
There is a bijection between tiles in the forward
partition and admissible paths.  The tile
corresponding to $a \to b$ is labelled by the
vertex pair $(v,w)$, where $v$ is the first
vertex of $a \to b$ and $w$ is the last one.
\end{lemma}

We let $T(a \to b)$ denote the tile in the
forward partition corresponding to the path
$a \to b$.  Let $\psi^*$ denote the pinwheel map,
defined in \S \ref{main result}.
 Here are our two main technical
results.

\begin{theorem}
\label{structure3}
Suppose $p \in T(a \to b)$ and
$q=\psi(p) \in T(c \to d)$ Then 
$(\psi^*)^k(q,b-1)=(q,c-1)$ for some $k \in \{0,...,n\}$.
\end{theorem}

\begin{theorem}
\label{structure2}
Let $a \to b$ be an admissible path,
labelled so that $b \geq a$. Then
$(\psi^*)^k(p,a-1)=(\psi(p),b-1)$
for some $k \in \{1,...,2n\}$.
\end{theorem}

In \S 5 and \S 6-7 we prove
Theorems \ref{structure3} and \ref{structure2}
respectively.  Combining these two results, we have

\begin{corollary}
\label{aux}
Let $T(a \to b)$ and $T(c \to d)$ be two
tiles such that there is a point $p$
such that $p \in T(a \to b)$ and
$\psi(p) \in T(c \to d)$.  Then there
is some integer $k \in \{0,...,3n\}$ such that
$(\psi^*)^k(p,a-1)=(\psi(p),c-1).$
\end{corollary}

Now we deduce the Pinwheel Theorem from
Corollary \ref{aux}.  The map $\psi$
is defined precisely on the interiors
of the tiles in the forward partition.
For $p$ in the interior of the tile
$T(a \to b)$, we define
\begin{equation}
\iota(p)=(p,a-1).
\end{equation}
As usual, we take the indices mod $n$.
The conclusion of Corollary \ref{aux}
is just a restatement of the equation
in the Pinwheel Theorem.

\subsection{Dynamical Corollaries}
\label{final}

\subsubsection{Proof of Corollary \ref{bijection}}

Let $p \in \R^2$ be a point with an
unbounded $\psi$-orbit $O$.
Let $p^* = \iota(p)$.  Let
$O^*$ be the $\psi^*$ orbit of $p^*$.
By the Pinwheel Theorem, there is an
infinite sequence $t_1,t_2,...$ such
that
$$\iota \circ \psi^k(p)=(\psi^*)^{t_k}(p^*).$$
This clearly shows that $O$ is unbounded if
and only if $O^*$ is unbounded.  
The analysis in Lemma \ref{far} shows that
$\pi(O^*-K^*) = O-K$ if $K$ is a sufficiently
large compact set.  Here $K^*=\pi^{-1}(K)$.

Since $\pi$ is a bijection outside of $K^*$, there
is only one orbit $O^*$ such that
$\pi(O^*-K)=O-K$.  Hence, the assignment $O \to O^*$,
which first seems to depend on the choice of $p$,
is well defined independent of the choice.
If $O_1$ and $O_2$ are different unbounded
orbits, they differ outside of $K$. Hence
$O_1^*$ and $O_2^*$ differ as well.
This shows that the assignment $O \to O^*$
is injective.  

Finally, let $O^*$ be some unbounded natural
orbit.  We just choose some $p \in O^*-K^*$ and
let $p=\pi(p^*)$.  Then the argument above
shows that $O^*$ is the image of $O$ under
our correspondence.  Hence, our correspondence
is a bijection.

\subsubsection{Proof of Corollary \ref{bijection2}}

In light of Corollary \ref{bijection}, we just
have to construct a bijection between the set of
forward unbounded natural orbits of $\psi^*$ and the
set of forward unbounded orbits of $\widehat \psi$.
But $\iota(\R^2)$ and $\widehat X$ agree outside
a compact set.  So, suppose that $O^*$ is a
forward unbounded natural orbit. The set
$\widehat O=O^* \cap \widehat X$
is a forward unbounded orbit of $\widehat X$.
The nature this construction
makes it clear that the correspondence
$O^* \to \widehat O$ is a bijection.

\subsubsection{Proof of Corollary \ref{bijection3}}

For each orbit $O$ of $\psi$, the intersection
$O \cap \Sigma_1$ is the corresponding orbit of
$\Psi$. In light of the analysis in Lemma \ref{far}
this gives a bijection between the set of unbounded
orbits of $\psi$ and the set of unbounded orbits
of $\Psi$. Similarly, there is a canonical bijection
between the set of unbounded orbits of $\widehat \psi$
and the set of unbounded orbits of $\widehat \Psi$.
Finally, Corollary \ref{bijection2} gives us a
canonical bijection between the unbounded orbits
of $\psi$ and the set of unbounded orbits of
$\widehat \psi$.  Composing all these bijections
gives us the desired result.

\newpage

\section{Spokes}

\subsection{Basic Definitions}
\label{leading}

Throughout the chapter, $P$ is a nice $n$-gon.
We say that a {\it minimal strip\/} is a strip
$\Sigma$ that contains the interior of $P$ in
its interior and is not a proper subset of any
other such strip. (The pinwheel strips are twice
as fat as certain minimal strips.)
We call an ordered pair $(v,w)$ of
vertices of $P$ a {\it maximal pair\/} if there is a 
minimal strip $\Sigma$ such that $v$ and $w$ lie in
distinct components of $\partial \Sigma$.
The spokes we defined in
\S \ref{spoke def} are precisely the line segments
joining vertices of maximal pairs.
As in the previous chapter, we let
$S_k$ denote the spoke corresponding to the vector $V_k$.

\begin{lemma}
\label{minimal1}
Let $S_1$ and $S_2$ be two consecutive spokes.  Then
$S_1$ and $S_2$ share a common vertex.  Moreover,
there is a minimal strip that contains all three
vertices of $\partial S_1 \cup \partial S_2$.
\end{lemma}

\startproof
For each $\theta \in \R \cup \infty$, there is a unique
minimal strip $\Sigma_{\theta}$ having slope $\theta$.
For all but $n$ value of $\theta$, the strip
$\Sigma_{\theta}$ contains $2$ vertices of $P$ on
its boundary.  For such value, there is a unique
spoke $S_{\theta}$ selected by $\Sigma_{\theta}$.
The assignment $\theta \to S_{\theta}$ is locally
constant and changes only at one of the $n$
special values.  At these special values,
$\Sigma_{\theta}$ contains $3$ vertices of $P$
in its boundary, and these are precisely the
endpoints of two consecutive spokes.  Our
lemma is a restatment of this geometric
picture.
\endproof

Recall that $e_1,...,e_n$ are the edges of $P$,
ordered according to slope.  The correspondence
$e_j \to S_j$ is what we call the
{\it clockwise correspondence\/}.
Figure 2.4 shows an example.

\begin{lemma}
\label{minimal2}
The clockwise correspondence is a bijection between edges and spokes.
\end{lemma}

\startproof
By the Pidgeonhole Principle, it suffices to prove that the
correspondence $e \to P_e$ is an injection.  So, for the sake
of contradiction, suppose that there are two edges
$e$ and $f$ such that $P_e=P_f=P$.  Let $a,b,c,d$ be
the edges of $e$ and $f$, ordered as in Figure 3.1.  
The vertex $b$ is farther from $\overline{cd}$ than is $a$
and the vertex $c$ is farther from $\overline{ab}$ than is
$d$.  This is a geometrically impossible situation.
\endproof

\begin{center}
\resizebox{!}{1.7in}{\includegraphics{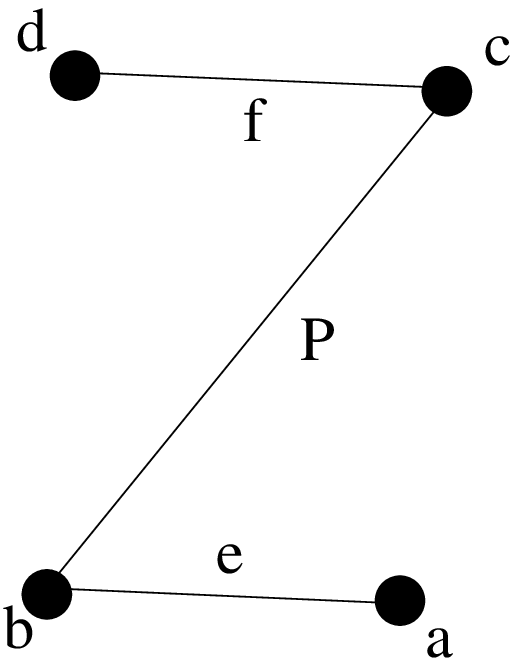}}
\newline
{\bf Figure 3.1:\/} An impossible situation
\end{center}

\subsection{Smooth Approximation}

We say that an {\it oval\/} is a smooth
and convex simple closed curve of strictly
positive  curvature. The condition
of strict positive curvature guarantees
that the tangent lines to the curve
vary strictly monotonically in small
neighborhoods.  This monotonicity is
all we really need.

Minimal strips and spokes are defined for
ovals just as we defined them for nice
polygons.
The spokes are easier to understand for
ovals than for polygons. Here are two
properties of spokes associated to
ovals.
\begin{enumerate}
\item Any oval has a unique spoke of a given slope.
\item Any two distinct spokes cross in their
interior.
\end{enumerate} 
These properties follow from
the monotonicity of the tangent lines
and continuity.  We omit the easy proofs.

  We are 
really only interested in the spokes of a nice
polygon, but one efficient way to understand
certain things about them is through smooth
approximations.
We can find a sequence $\{P_n\}$ of
ovals that converges to $P$ in the Hausdorff topology.
More precisely, we mean that $P_n$ is contained
in the $(1/n)$-tubular neighborhood of $P$, and
{\it vice versa\/}.   We omit the straightforward
proof of this fact.

\begin{lemma}
For any spoke $S$ of $P$ there is a sequence
$\{S_n\}$ that converges to $S$.  Here
$S_n$ is a spoke of $P_n$.
\end{lemma}

\startproof
There is a spoke $S_n$ of $P_n$ that has
the same slope as $S$. Evidently $S_n$
converges to $S$.
\endproof

\subsection{Admissible Paths and Vertex Pairs}
\label{define admissible}

Here we give a quick illustration of the utility of
smoothly approximating nice polygons by ovals.
(We will give somewhat more substantial applications
in later sections of this chapter.)

\begin{lemma}
\label{spoke1}
Any two spokes of a nice polygon $P$ intersect.
\end{lemma}

\startproof
The corresponding result for ovals is fairly
obvious.  We mentioned it above.  The corresponding
result for nice polygons follows by continuity.
\endproof

We say that a
pair $(v,x)$ of vertices of $P$ is
{\it admissible\/} if there is some vertex
$w$ such that $(v,w)$ is a maximal pair and
The clockwise path joining $v$ to $w$ contains $x$.
Figure 3.2 shows a picture.

\begin{center}
\resizebox{!}{2in}{\includegraphics{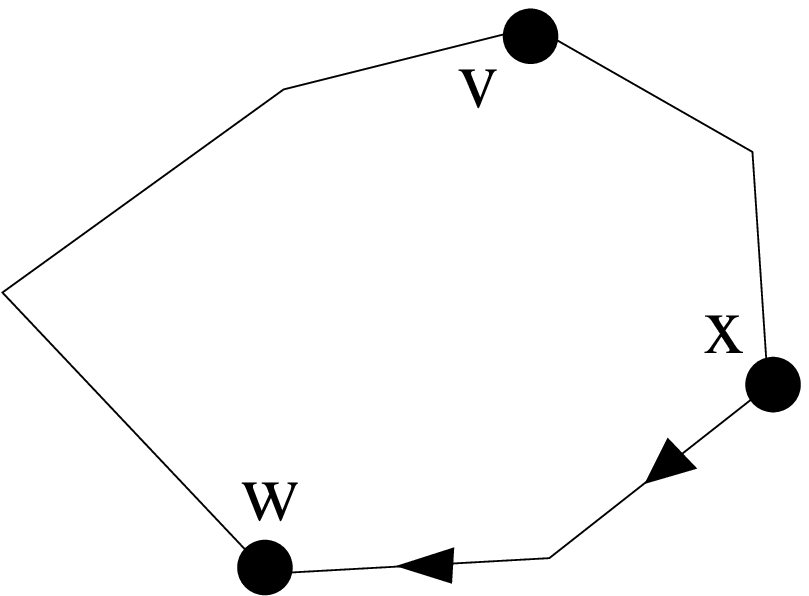}}
\newline
{\bf Figure 3.2:\/} Admissible pairs of vertices.
\end{center}

\begin{lemma}
Suppose that $(v,w)$ and $(w,v)$ are both
admissible pairs.  Then $(v,w)$ is a maximal pair.
\end{lemma}

\startproof
If $(v,w)$ is an admissible pair that is not maximal,
then there are maximal pairs $(v,w')$ and $(w,v')$
that relate to $(v,w)$ and $(w,v)$ as described in
the definition of admissible pairs. By construction,
these two maximal pairs do not intersect.  This
contradicts.  Lemma \ref{spoke1}.
\endproof

The next lemma refers to admissible paths,
which we defined in \S \ref{spoke def}.

\begin{lemma}
\label{connect}
A pair of vertices $(v,w)$ is admissible if and only
if $v$ and $w$ respectively are the starting and
endpoint points of an admissible path.
\end{lemma}

\startproof
Each maximal pair of vertices is clearly admissible.
These correspond to single spokes.  Conversely, each
admissible path of length $1$ is an oriented spoke
and hence corresponds to a maximal pair of vertices.
If we start with 
an admissible path $a \to b$ and minimally lengthen
it to the new admissible path $a \to b'$, the new
endpoint is a vertex adjacent to, and counterclockwise
from, the old endpoint.   At the same time, the
admissible vertices are obtained from the maximal
vertices by moving the endpoints counterclockwise.
Our result follows from these facts and from induction.
\endproof

\subsection{Interpolation Properties of Spokes}
\label{interpolation}
\label{harmony}

Suppose that $S_a$ and $S_b$ are two spokes
of $P$.  The indices correspond to the
ordering on the spokes that comes from
their slopes.  We take indices mod $n$
and arrange that $a<b$.  Let
$S_a(1)$ and $S_a(2)$ be the two vertices
of $S_a$.  Let $S_b(1)$ and $S_b(2)$ be
the two vertices of $S_b$.  We choose
these labels so the vertices
$$S_a(1); \hskip 20 pt
S_b(1); \hskip 20 pt S_a(2); 
\hskip 20 pt S_b(2)$$
are counterclockwise cyclically ordered, as in
Figure 3.3.

\begin{center}
\resizebox{!}{1.8in}{\includegraphics{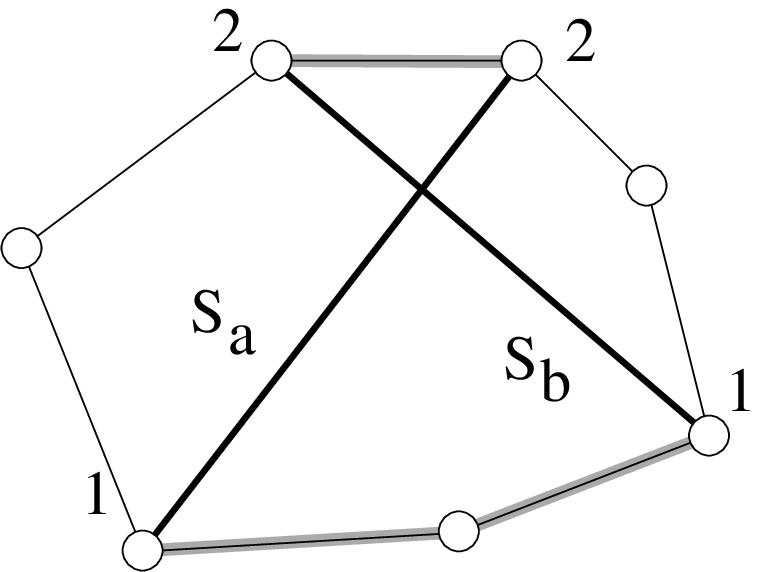}}
\newline
{\bf Figure 3.3:\/} Two spokes and two arcs.
\end{center}

We distinguish $2$ arcs of $P$. The first one, $V(a,b,1)$,
connects $S_a(1)$ to $S_b(1)$ while avoiding
$S_a(2)$ and $S_b(2)$.  We define $V(a,b,2)$
the same way, with the roles of $(1)$ and $(2)$
interchanged.   These two paths are highlighted
in grey in Figure 3.3.  The way we have assigned
spokes to edges in the previous section is not
symmetric with respect to order reversal.
This subtle asymmetry makes itself felt in the
``lopsided'' condition $a \leq j<b$ in the next result.

\begin{lemma}
\label{between}
Suppose that $a \leq j< b$. Then
$e_j \subset V(a,b,1) \cup V(a,b,2)$.
Here $e_j$ is the edge associated to
the spoke $S_j$.
\end{lemma}

\startproof
The same kind of limiting argument as in Lemma \ref{spoke1}
shows that the endpoints of $S_j$ lie in
$V(a,b,1) \cup V(a,b,2)$ for each relevant $j$.
This lemma fails only if 
an endpoint of $S_j$ is also an 
endpoint of $V(a,b,1) \cup V(a,b,2)$ and the
edge $e_j$ sticks out the side.  Figure 3.4 shows
one of the two possibilities.

\begin{center}
\resizebox{!}{3in}{\includegraphics{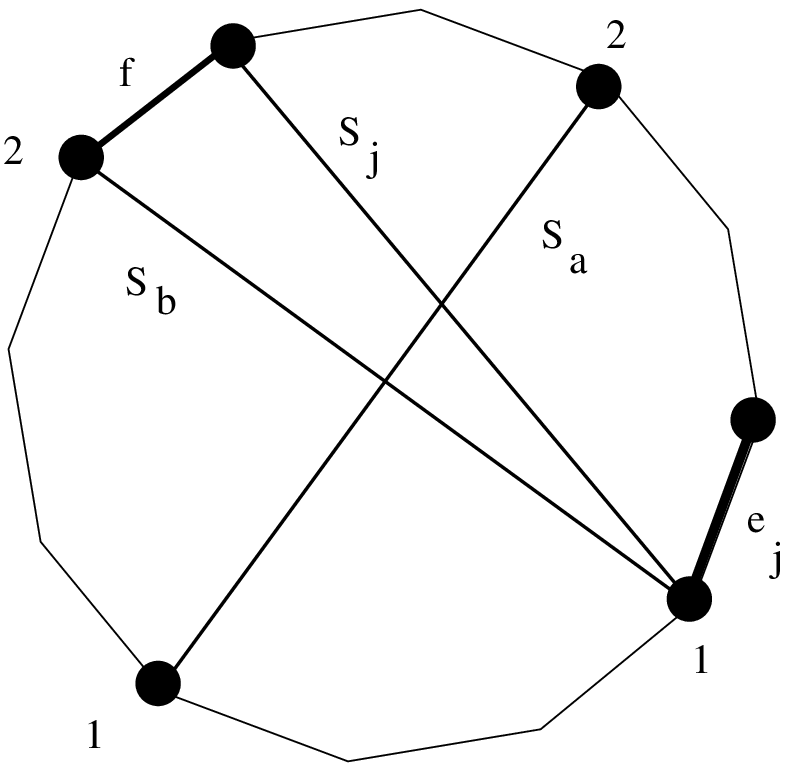}}
\newline
{\bf Figure 3.4:\/} The edge sticking out the side.
\end{center}

As $j$ increases to $b$, the spoke $S_j$ turns counterclockwise,
but its endpoints cannot leave $V(a,b,1) \cup V(a,b,2)$.  One
of the endpoints is already stuck at the right endpoint
of $V(a,b,1)$ and so this endpoint cannot move at all.
The only possibility is that the spokes $S_{j+1},...,S_b$
all share the endpoint $S_b(1)$.  In particular, we may
assume that $j=b-1$.  In this case, the two unequal endpoints
of $S_j$ and $S_b$ are the endpoints of an edge, say $f$.
Since $S_j$ and $S_b$ are consecutive, it follows from
Lemma \ref{minimal1} that there is a minimal strip whose
boundary contains $\partial S_j \cup \partial S_b$.  One
of the sides of this strip must contain $f$. But then,
the clockwise correspondence assigns $S_j$ to $f$.  On the
other hand, the clockwise correspondence assigns $S_j$ to $e_j$
and $e_j \not = f$.
This contradicts Lemma \ref{minimal2}, which says that
the clockwise correspondence is a bijection.
\endproof

Let $p_1,p_2,p_3 \in \R^2-P$ be a portion of an
outer billiards orbit.   We say that the triple
$(p_2,S_a,S_b)$ is {\it harmonious\/} if the
ray $\overrightarrow{p_1p_2}$ is tangent to
$P$ at an endpoint of $S_a$ and the ray
$\overrightarrow{p_2p_3}$ is tangent to
$P$ at an endpoint of $S_b$.  Figure 3.5 shows
the situation.

\begin{center}
\resizebox{!}{3.5in}{\includegraphics{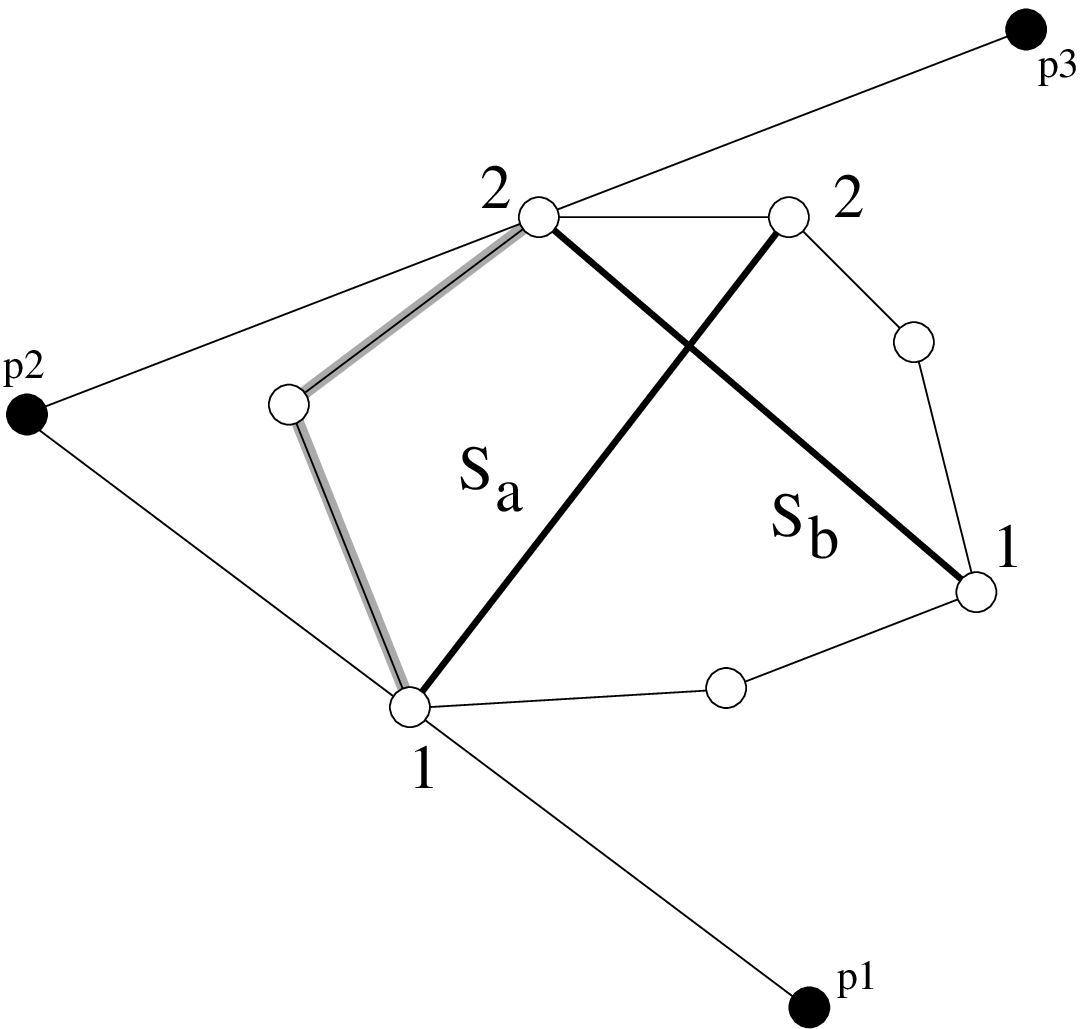}}
\newline
{\bf Figure 3.5:\/} A harmonious triple.
\end{center}

Figure 3.5 also shows a distinguished arc $A(p)$ having
endpoints $S_a(1)$ and $S_b(2)$.  We say that this
$A(p)$ is the arc {\it subtended\/} $p$.  The arc
$A(p)$ consists of those points $q \in P$ such that
the line
segment joining $q$ to $p$ only intersects $P$ at $q$.
Say that two closed arcs of $P$ are
{\it almost disjoint\/} if they share an endpoint.

\begin{lemma}
\label{subtend}
Suppose $(p,S_a,S_b)$ are a harmonious triple.
Then $A(p)$ is almost disjoint from $V(a,b,1)$
and $V(a,b,2)$.
\end{lemma}

\startproof
All the same definitions make sense for
ovals, and the result there is obvious.
The result now follows for nice polygons
by continuity.
\endproof

\subsection{The Pinwheel Orientation}
\label{pinwheel orientation}

We define a canonical orientation on the
spokes of a nice $n$-gon $P$.
Each spoke $S_j$ corresponds to a unique pinwheel
strip $\Sigma_j$ in the sense that the
vector $V_j$ associated to $\Sigma_j$ points
from one endpoint of $S_j$ to the other.
We orient $S_j$ by saying that ${\rm oriented\/}(S_j)=V_j$,

We say that a spoke $S_j$ is {\it special\/}
if the three spokes $S_{j-1}, S_j, S_{j+1}$ share a common vertex,
and otherwise {\it ordinary\/}.  For instance,
in Figure 2.2 the spokes $S_1$ and $S_3$ are special and the rest are
ordinary. 

\begin{lemma}
\label{orientation}
Let $\gamma$ be an admissible path of length
at least $3$.  The orientation on the 
spokes of $\gamma$ induced by the orientation on
$\gamma$ coincides with the pinwheel orientation
on all edges but the last one.  On the last
spoke, the two orientations agree if and only if
the spoke is ordinary.
\end{lemma}

\startproof
Consider a subpath of $\gamma$ that contains two spokes,
say, $S_1$ and $S_b$.  Here $b \in \{2,3,4...\}$.
Figure 3.6 shows the possible pictures, depending on $b$.

\begin{center}
\resizebox{!}{1.5in}{\includegraphics{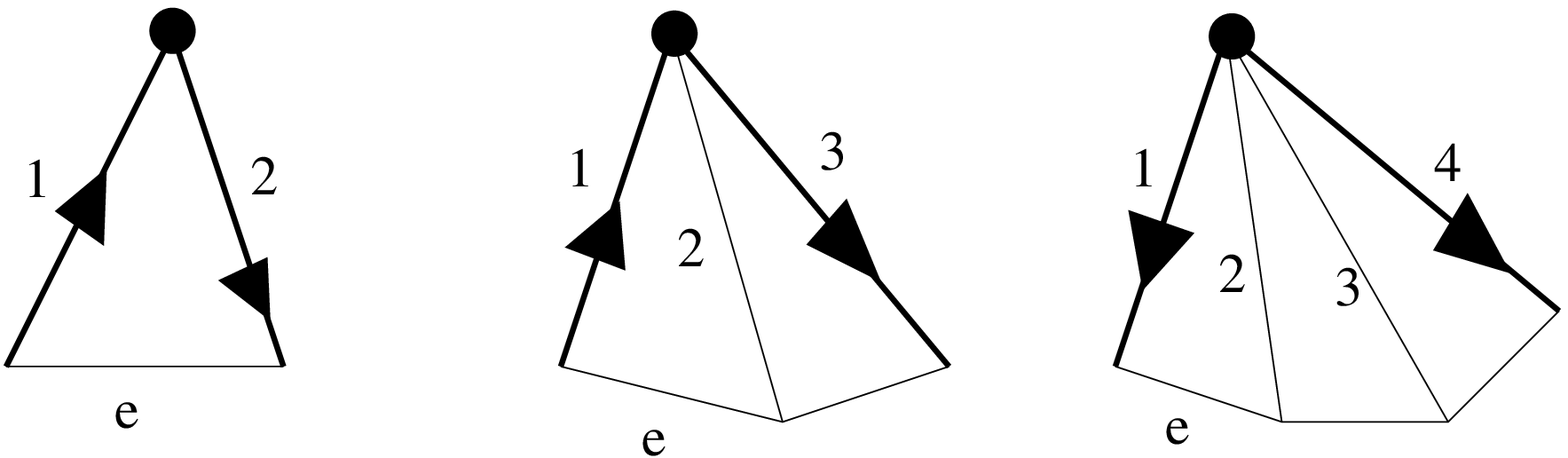}}
\newline
{\bf Figure 3.6:\/} The possibilities for $b=2,3,4$.
\end{center}

In all cases, it follows from Lemma \ref{minimal1} that
there exists a minimal strip $\Sigma'$ such that $\partial \Sigma'$
contains $\partial S_1 \cup \partial S_2$.  In particular,
$\partial \Sigma'$ contains the edge $e$ shown in Figure 3.6.
But then the strip $\Sigma$ that is twice as wide as $\Sigma'$
is a pinwheel strip containing $e$ in its boundary.
Comparing Figures 1.2 and 3.6, we see that the orientation
given to $S_1$ by $\gamma$ coincides with the pinwheel
orientation on $S_1$.

It remains to describe what happens for the last spoke of $\gamma$.
The argument we just gave never uses the third property
satisfied by an admissible path.  In case the last spoke of
$\gamma$ is ordinary, we can prolong $\gamma$ by two more
spokes in such a way that the longer path $\gamma'$ satisfies
the first two properties for an admissible path.  The 
previous argument then shows that the pinwheel orientation
and the $\gamma$-induced orientations coincide the third-to-last
spoke of $\gamma'$, which is the last spoke of $\gamma$.

\begin{center}
\resizebox{!}{1.6in}{\includegraphics{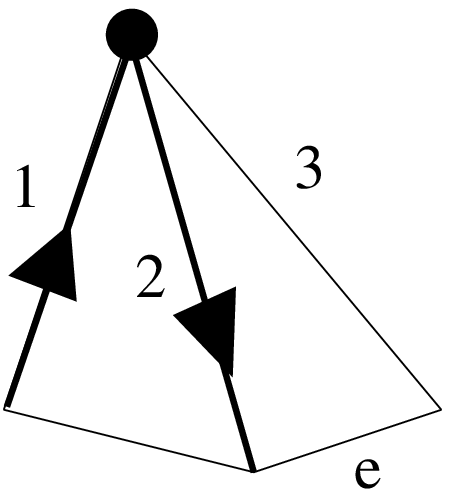}}
\newline
{\bf Figure 3.7:\/} Ending on a special spoke.
\end{center}

Suppose that the last spoke of $\gamma$ is special.
Just for the sake of argument, let's suppose that
the last spoke of $\gamma$ is labelled $S_2$.
Figure 3.7 shows the picture of the last two
spokes of $\gamma$ when $S_2$ is special.
In this case $e \in \Sigma_2$.  Comparing
Figures 1.2 and 3.7 we see this time that the
orientation on $S_2$ given by $\gamma$ is
$-V_j$, as claimed.
\endproof

We also record the following result.
\begin{lemma}
\label{special}
The edges $e_0$ and $e_1$ are adjacent if and only
if $S_1$ is a special spoke.
\end{lemma}

\startproof
If $S_1$ is a special spoke, then the spokes
$S_0,S_1,S_2$ all share a common point. The
other three points are the vertices of the
arc $e_0 \cup e_1$.  Hence $e_0$ and $e_1$ are
adjacent.  Conversely, if $e_0$ and $e_1$ are
adjacent, then it follows from Lemma \ref{minimal1}
that $e_0$ is the edge connecting two
points of $\partial S_0 \cup \partial S_1$ and
$e_1$ is the edge connecting two
points of $S_1 \cup S_2$.   But then
$S_0,S_1,S_2$ share a common vertex.
Hence $S_1$ is special.
\endproof

\newpage

\section{The Forward Partition}

\subsection{The Combinatorics of the Partition}
\label{tile geometry}

We say that a {\it forward wall\/} is a line segment or ray
that lies in the closure of $2$ tiles of the forward
partition.  We say that a forward wall $W$ is {\it primary\/}
if the outer billiards map is not defined on interior
points of $W$.  We say that $W$ is {\it secondary\/}
if $W$ is not primary. For a point in the interior
of a secondary wall, the outer billiards map is
defined but the next iterate is not defined.
Put another way, the secondary walls are obtained
by pulling back the union of the primary walls under
the outer billiards map.
Figure 4.1 shows pattern for a pentagon.
The primary walls divide  $\R^2-P$ into $n$ distinct
{\it primary cones\/}.

\begin{center}
\resizebox{!}{1.8in}{\includegraphics{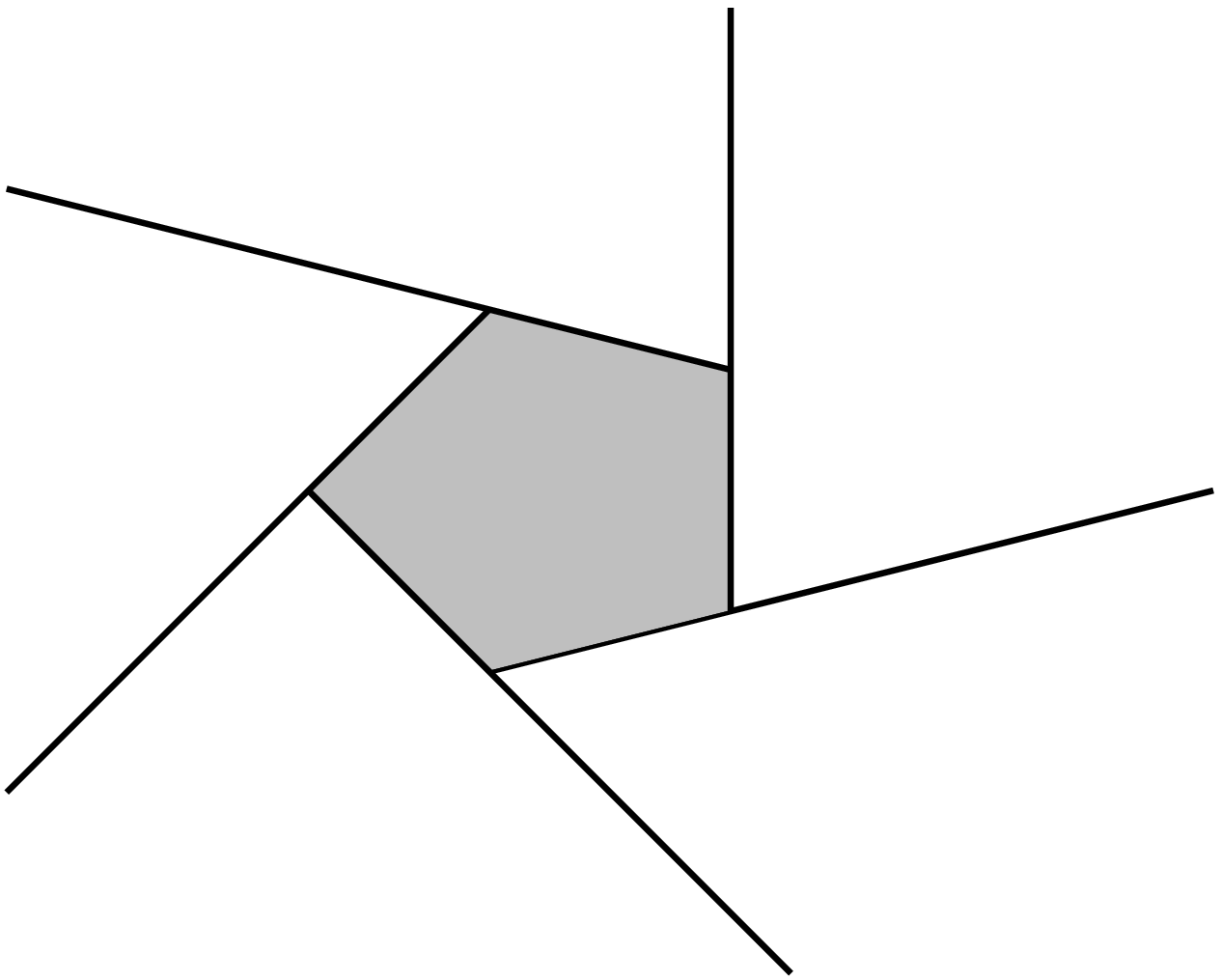}}
\newline
{\bf Figure 4.1:\/} Primary walls
\end{center}

Each secondary wall is either a line segment
or a ray. 

\begin{lemma}
An endpoint of a secondary wall liea on a primary call and
cannot be the vertex of the cone that contains it.
\end{lemma}

\startproof
Let $p_0$ be the endpoint of a secondary wall $w_0$.
Let $p_1$ be the image of $p_0$ under the outer billiards map.
If $p_0$ is in the interior of a cone then
$p_1$ must lie on an interior point of
a primary wall $w_1$.  Pulling back a neighborhood
of $p_1$ in $w_1$, we see that $w_0$ contains a
$2$-sided neighborhood of $p_0$.  Hence
$p_0$ is not an endpoint of $w_0$.
The second statement follows from the fact that
the square-outer billiards map is clearly
defined in a neighborhood of the vertex of
a cone.
\endproof

\begin{lemma}
No two secondary walls intersect.
\end{lemma}

\startproof
Suppose that $p_0 \in w_1 \cap w_2$ is a point of intersection
between two secondary walls. We can approximate
$p_0$ by a sequence $\{p_0(k,n)\}$ of points in $w_k$
for $k=1,2$.  We can arrange that both sequences
lie in the interior of the same cone.  Hence, the
outer billiards map is defined on both sequences.
Let $p_1(k,n)$ denote the image of $p_0(k,n)$
under the outer billiards map.  Since the outer
billiards map is an isometry on the interior of
each cone, the two sequences $\{p_1(k,n)\}$
converge to the same point $q_1$.  By
construction, $q_1$ lies on two different primary
walls.  But the primary walls only intersect at the
vertices of $P$. Hence $q_1$ is a vertex of $P$.
Given the definition of the outer billiards map,
$p_0$ must be the same vertex of $P$.  This is
ruled out by the previous lemma.
\endproof

According to the lemmas above, the secondary walls
cut across the cones like crosshatching, as
shown schematically in Figure 4.2.  
Every bounded tile is either a triangle or a quadrilateral.
The triangular bounded tiles contain a neighborhood
of a cone apex.  The quadrilateral bounded tiles
have one pair of opposite sides that are primary
walls and another pair of opposite sides that are
secondary walls. As we justify in the proof of Lemma \ref{maxx} below,
each secondary wall is parallel to an edge of $P$,
and the pattern is as indicated by the labels in
Figure 4.2.  Compare Figure 2.2.

\begin{center}
\resizebox{!}{2.6in}{\includegraphics{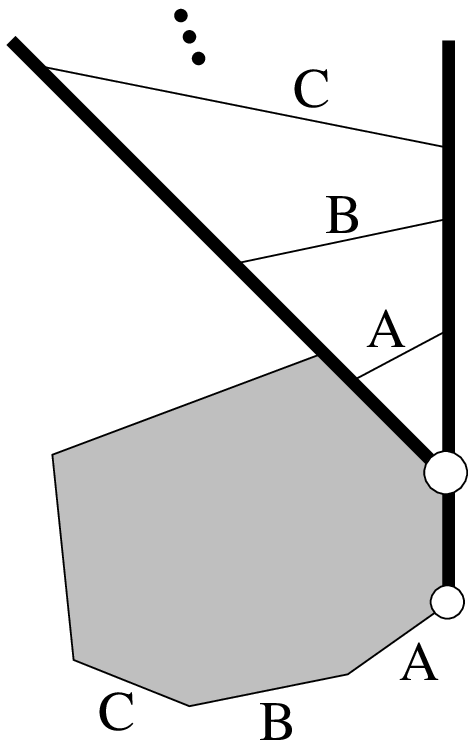}}
\newline
{\bf Figure 4.2:\/} Secondary walls.
\end{center}

\subsection{Proof of Lemma \ref{structure1}}
\label{admissible proof}

Just for this section we say that the pair
$(v,w)$ of vertices is a {\it label\/} if
it labels some nonempty tile in the forward partition.

\begin{lemma}
\label{maxx}
Every maximal pair of vertices is a label.
\end{lemma}

\startproof
Let $(v,w)$ be a maximal pair, and let
$\Sigma$ be the corresponding pinwheel strip.
We rotate so that $\Sigma$ is
horizontal, as in Figures 4.3 and 4.4 below.
If $p \in \R^2$ lies just below $L$ and far to the right, then 

\begin{equation}
\label{local1}
\psi(p)=p+V.
\end{equation}

\begin{center}
\resizebox{!}{1.3in}{\includegraphics{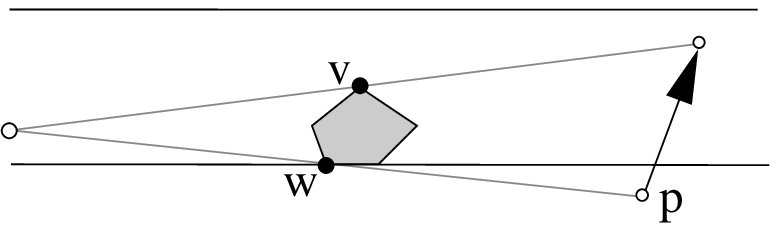}}
\newline
{\bf Figure 4.3:\/} $(w,v)$ is a label.
\end{center}

If $p$ lies just above $L'$ and far to the
left, the map $\psi$ has the form 
\begin{equation}
\label{local2}
\psi(p)=p-V.
\end{equation}  

\begin{center}
\resizebox{!}{1.3in}{\includegraphics{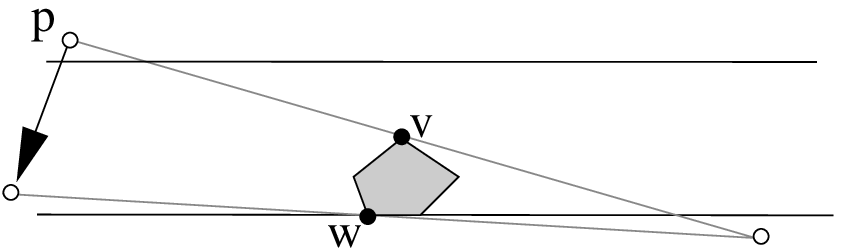}}
\newline
{\bf Figure 4.4:\/} $(v,w)$ is a label.
\end{center}

It is a consequence
of Equations \ref{local1} and \ref{local2} that
there is an unbounded region labelled by
$(w,v)$ and an unbounded region labelled by $(v,w)$.
\endproof

Now we show that the set of labels coincides with the
set of admissible pairs.
Let's first show that such any admissible pair $(v,x)$ is a label.
First choose a point $p_0 \in \R^2-P$ in the unbounded tile
labelled by $(v,w)$.   Consider the outer billiards
orbit $p_0,p_1,p_2$, with $p_2=\psi(p_0)$ as usual.
The situation is shown in Figure 4.4.

\begin{center}
\resizebox{!}{3.3in}{\includegraphics{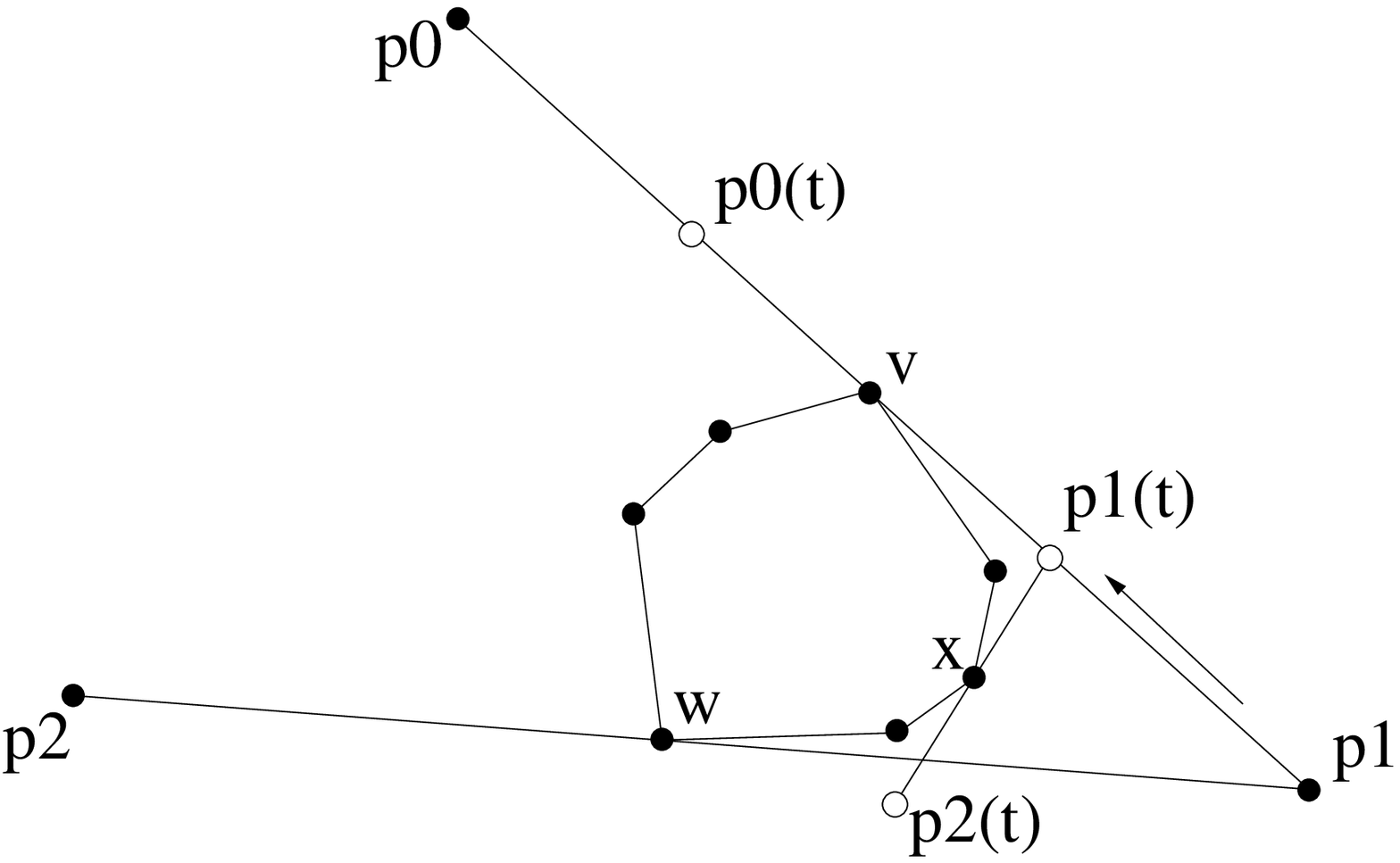}}
\newline
{\bf Figure 4.5:\/} $(v,x)$ is a label.
\end{center}

Let 
$p_1(t)$ denote a point that varies along the line
segment connecting $p_1$ to $v$.  Let
$p_0(t)$ and $p_2(t)$ be the backwards and
forwards images of $p_1(t)$ under the outer billiards map.
As $p_1(t)$ varies from $p_1$ to $v$, the moving
tangent line from $p_1(t)$ to $P$ encounters every
vertex along the oriented path from $v$ to $w$.  Hence, there is some
value of $t$ for which $p_0(t)$ lies in a region
labelled by $(v,x)$.  Hence $(v,x)$ is a label.
This argument, incidentally, justifies the comments
about the secondary walls mentioned at the end
of the last section.

To show that all the labels are admissible, we
observe that we can reverse the process
we just described, starting with a point $p_0$
in the region labelled by $(v,x)$ and then
moving outward along the ray pointing from $v$ to $p_0$
until we arrive at a point that lies in a region
labelled by the maximal pair $(v,w)$.  This
shows that $(v,x)$ is admissible.

Combining the result here with Lemma \ref{connect},
we establish Lemma \ref{structure1}.

\subsection{Some Geometric Estimates}

The bounded tiles in a single cone $C$ all have a
common boundary with one of the strips $\Sigma$, but
are disjoint from that strip.  If we rotate so
that $\Sigma$ is horizontal, as in Figure 4.6, then
the shaded region $T$ beneath $\Sigma$ is the union
of the bounded tiles in $C$.

\begin{center}
\resizebox{!}{3.4in}{\includegraphics{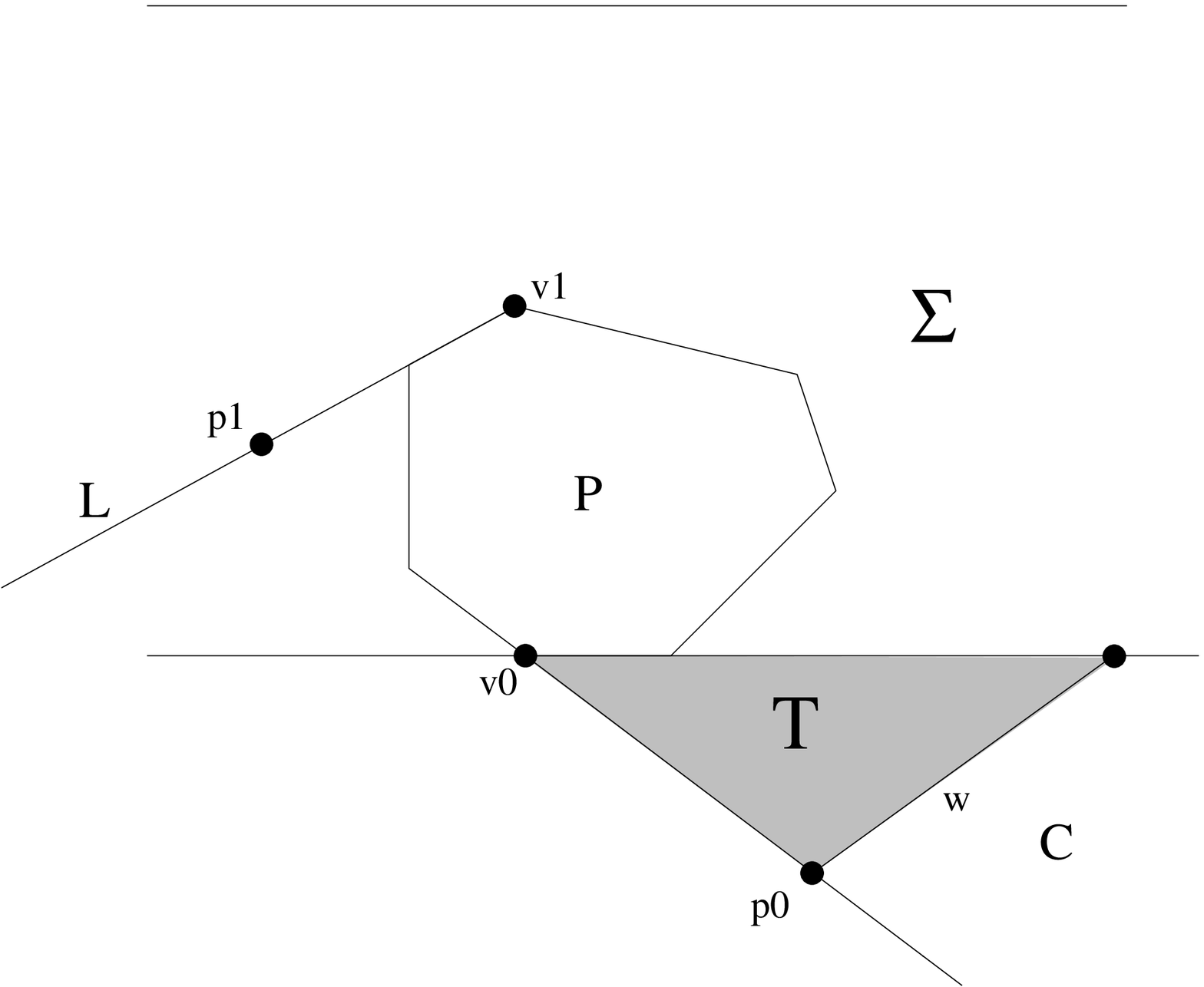}}
\newline
{\bf Figure 4.6:\/} Bounding the tiles.
\end{center}

\begin{lemma}
\label{small}
No point of $T$ is farther from $\partial \Sigma$ than
half the width of
$\Sigma$.
\end{lemma}

\startproof
The boundary wall $w$ of $T$ opposite the cone apex $v_0$
is the largest bounded
secondary wall in $C$.  The secondary wall $w$ 
is parallel to a line $L$ which
extends one of the sides of $P$.  The
fact that $w$ is bounded gives control on the
slope of $L$, and the basic geometric fact is that
the point $p_1$ necessarily lies in the lower half
of $\Sigma$.  On the other hand, $p_0$ and $p_1$
are equidistant from the lower boundary of $\Sigma$.
Here $p_1$ is the image of $p_0$ under the outer
billiards map.
But this shows that the distance from $p_0$ to
the lower boundary of $\Sigma$ is less than half
the width of $\Sigma$.  Since $p_0$ is the
extreme point of $T$, the same statement can be
made for any point of $T$.
\endproof

Let $e_1,...,e_n$ be the edges of $P$.
Let $\partial_1\Sigma_j$ be the component of
$\partial \Sigma_j$ that contains $e_j$. Let
$\partial_2 \Sigma_j$ be the other component.

\begin{lemma}
\label{sep}
If $e_{0}$ is adjacent to $e_1$ then $T(1 \to b)$ and
$P$ lie on the same side of $\partial_1 \Sigma_{0}$.
Otherwise, $T(1 \to b)$ and $P$ lie on the same side of 
$\partial_2 \Sigma_0$.
\end{lemma}

\startproof
If $e_0$ and $e_1$ are adjacent, then $e_0=f_0$ in Figure 4.7.
The lines extending
$e_0$ and $e_1$ contain the boundary of the primary cone
that contains $T(a \to b)$.  This is the shaded
region in Figure 4.7.  In this case, the result is
obvious.  In the second case, $e_0=g_0$,
and we want to rule out the possibility
that $T(a \to b)$ intersects the black cone.
But one checks easily that
points lie in the black cone lie in tiles labelled by
vertex pairs $(x,y)$
where $y$ lies in the interior of the
clockwise arc connecting the two endpoints of the
spoke $S_a$.  This is a contradiction.
\endproof

\begin{center}
\resizebox{!}{3.9in}{\includegraphics{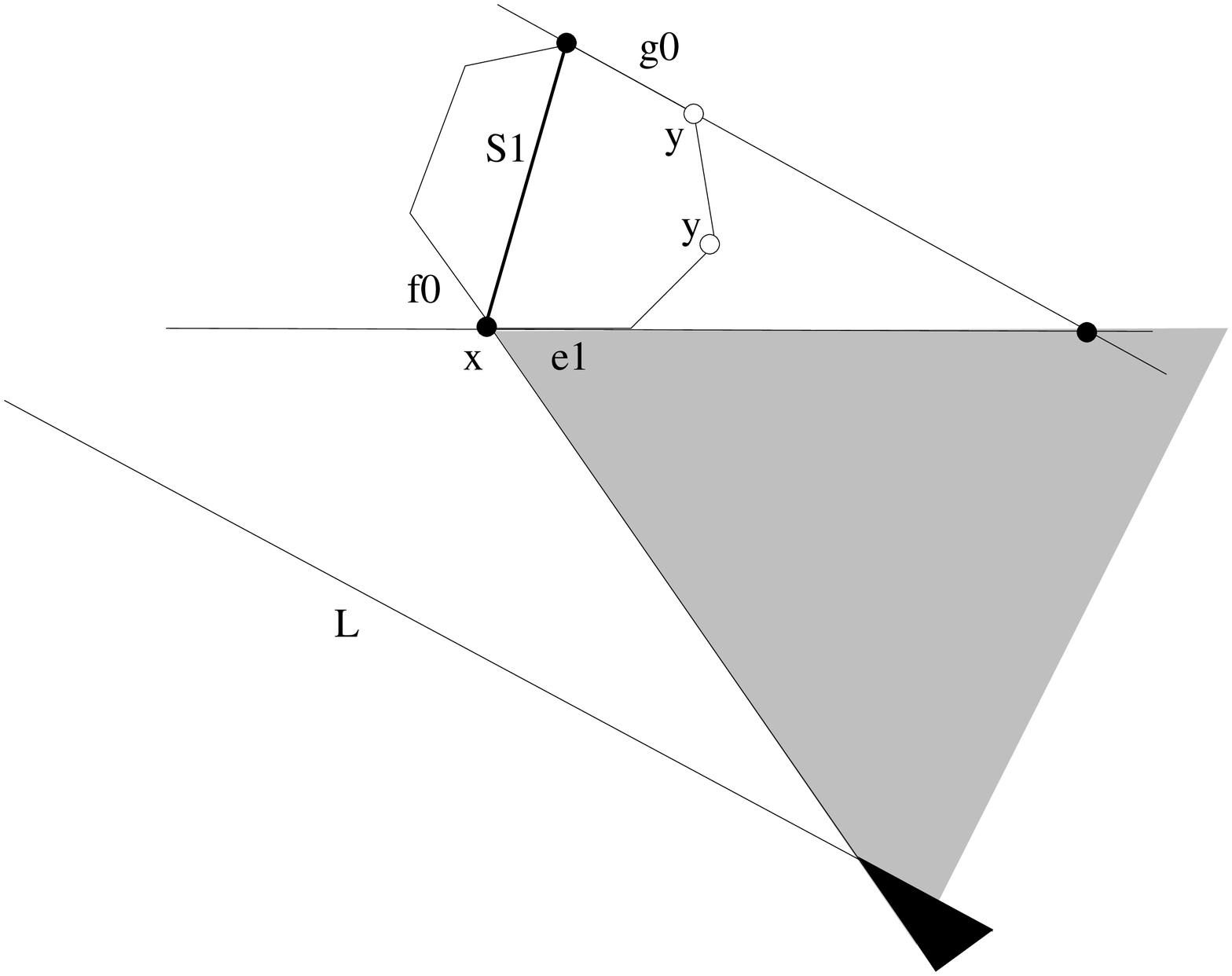}}
\newline
{\bf Figure 4.7:\/} The excluded region.
\end{center}

Recall that the vectors $V_1,...,V_n$ are 
associated to the pinwheel strips $\Sigma_1,...,\Sigma_n$.
In our proof of Theorem \ref{structure2} we will
be interested in sets such as
\begin{equation}
T(a \to b;a) = T(a \to b)+2V_a.
\end{equation}

\begin{corollary}
\label{prepin2}
If $e_{0}$ is adjacent to $e_1$ then
the line $\partial_2\Sigma_{0}$ separates $P$ from $T(1 \to b;1)$.
Otherwise, the line $\partial_1 \Sigma_0$ separates $P$ from $T(1 \to b;1)$.
\end{corollary}

\startproof
Suppose that $e_0$ and $e_1$ are adjacent.
Since $e_0$ and $e_1$ are adjacent, the vector $V_1$ joins a
a point on $\partial_1 \Sigma_0$ to a point on the
centerline of $\Sigma_0$.  Hence
$$\partial_1 \Sigma_0+2V_1 = \partial_2 \Sigma_0.$$
The first case of this lemma now follows from 
Lemma \ref{sep}.  The point is that adding $V_1$
ejects $T(1 \to b)$ outside of $\Sigma_0$,
and onto the correct side.

Suppose that $e_0$ and $e_1$ are not adjacent.
This time $V_1$ joins a point on the centerline of
$\Sigma_0$ to a point on $\partial_2 \Sigma_0$.
Hence $\partial_2 \Sigma_0+2V_1=\partial_1\Sigma_0$.
The rest of the proof is the same in this case.
\endproof

Recall from \S \ref{orientation} that a spoke $S_a$ is
{\it special\/} if the spokes $S_{a-1},S_a,S_{a+1}$
share a common vertex.  Otherwise, $S_a$ is ordinary.
Combining Lemma \ref{prepin2} with Lemma \ref{special} and
relabelling, we have the following result.

\begin{corollary}
\label{prepin3}
If $S_a$ is a special spoke, then
the line $\partial_2\Sigma_{a-1}$ separates $P$ from $T(a \to b;a)$.
Otherwise, the line
$\partial_1 \Sigma_{a-1}$ separates $P$ from $T(a \to b;a)$.
\end{corollary}

The last result we prove is not needed
anywhere in this paper. We include it because we
think gives a nice piece of extra information
about the geometry of the partition.

\begin{lemma}
\label{exit}
A partition tile $T$ is unbounded if
and only if $\psi(T) \cap T \not = \emptyset$.
\end{lemma}

\startproof
If $T$ is unbounded, then $T$ contains an infinite cone.
Being a uniformly bounded piecewise translation, $\psi$
cannot map an infinite cone off itself.
It remains to show that $\psi(T) \cap T \not = \emptyset$
implies that $T$ is unbounded.
Let $(v,w)$ be the pair of
vertices labelling $T$. It suffices to prove
that $(v,w)$ is a maximal pair.   Let
$p_0,p_1,p_2...$ be the forwards outer billiards
orbit of a point $p_0 \in T$,  Suppose
$p_2=\psi(p_2) \in T$.

\begin{center}
\resizebox{!}{3in}{\includegraphics{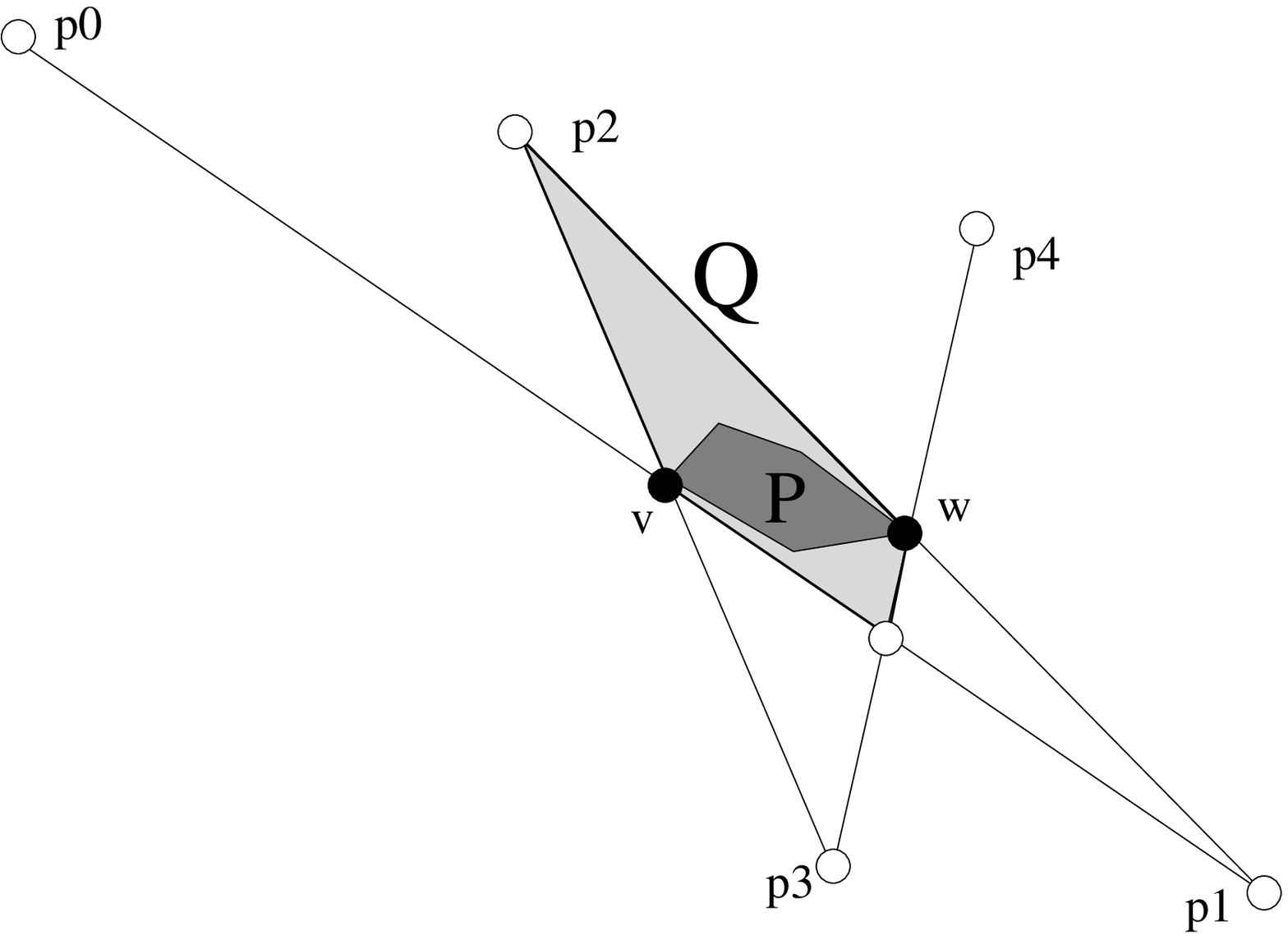}}
\newline
{\bf Figure 4.8:\/} A special quadrilateral
\end{center}

By construction, $P \subset Q$, where $Q$ is
the quadrilateral with vertices
$v$, $w$, $p_2$, and
$\overline{p_0p_1} \cap \overline{p_3p_4}$.
Since $v$ and $w$ are opposite vertices of $Q$, 
there is an infinite strip that contains
$v$ and $w$ on its boundary and containes the
rest of $Q$ in its interior.  This strip picks
out $(v,w)$ as a maximal pair.
\endproof

\newpage

\section{Proof of Theorem \ref{structure3}}
\label{proof3}

\subsection{Reformulation of the Result}

We can define the backwards partition
for $P$ just as we defined the forwards
partition.  We just use the inverse map
$\psi^{-1}$.  The map $\psi$ sets up
a bijection between the tiles in the forward
partition and the tiles in the backward
partition.  A tile in the forward partition
labelled by a pair of vertices $(v,w)$
corresponds to a tile in the backward
partition labelled by a pair of vertices $(w,v)$.

We change our notation so that
$T_+(a \to b)=T(a \to b)$, and
$T_-(a \to b)$ denotes a tile in the backwards
partition corresponding to the ``backwards admissible path''
The backwards admissible paths have the same definition
as the forwards admissible paths, except that the
spokes are traced out in reverse cyclic order.
In short, if we reverse a forwards admissible
path, we get a backwards admissible path.
From this structure, we have

\begin{equation}
\label{reverse}
\psi(T_+(a \to b)) = T_-(b \to a).
\end{equation}
We define
\begin{equation}
T_-(b) = \bigcup_a T_-(b \to a); \hskip 40 pt
T_+(c) = \bigcup_d T_+(c \to d).
\end{equation}

If $c=b$ there is nothing to prove.
If necessary, we add $n$ to $c$ so that
$b<c$.  In this chapter, we will prove the
following result.
\begin{equation}
\label{goal2}
T_-(b) \cap T_+(c) \subset \Sigma_j; \hskip 20 pt j=b,...,(c-1).
\end{equation}
Theorem \ref{structure3} is a quick corollary.
\newline
\newline
\noindent
{\bf Proof of Theorem \ref{structure3}:\/}
Our point $q$ in Theorem \ref{structure3} lies $T_-(b) \cap T_+(c)$.
Hence $$q \in \Sigma_j; \hskip 30 pt j=b,...,(c-1).$$
Hence
$$\psi^*(q,j)=(q,j+1); \hskip 30 pt j=(b-1),...,(c-2).$$  Hence
$(\psi^*)^k(q,b-1)=(q,c-1)$ for $k=b-a<n$.
\endproof

\subsection{A Result about Strips}

Let $P$ be a nice polygon. Let $\Sigma$ the pinwheel strip
associated to an edge $e$ of $P$.  We rotate and scale
so that $\Sigma$ is bounded by the lines $y=0$ and $y=1$,
as in Figure 5.2.  Let $v$ be the vertex
of $P$ farthest from the bottom edge $e$ of $P$.  We translate
so that $v=(0,1/2)$. Let
$T$ be the triangle formed by lines extending the
two edges of $P$ incident to $v$ and the top boundary
component of $\Sigma$. 

\begin{center}
\resizebox{!}{1.8in}{\includegraphics{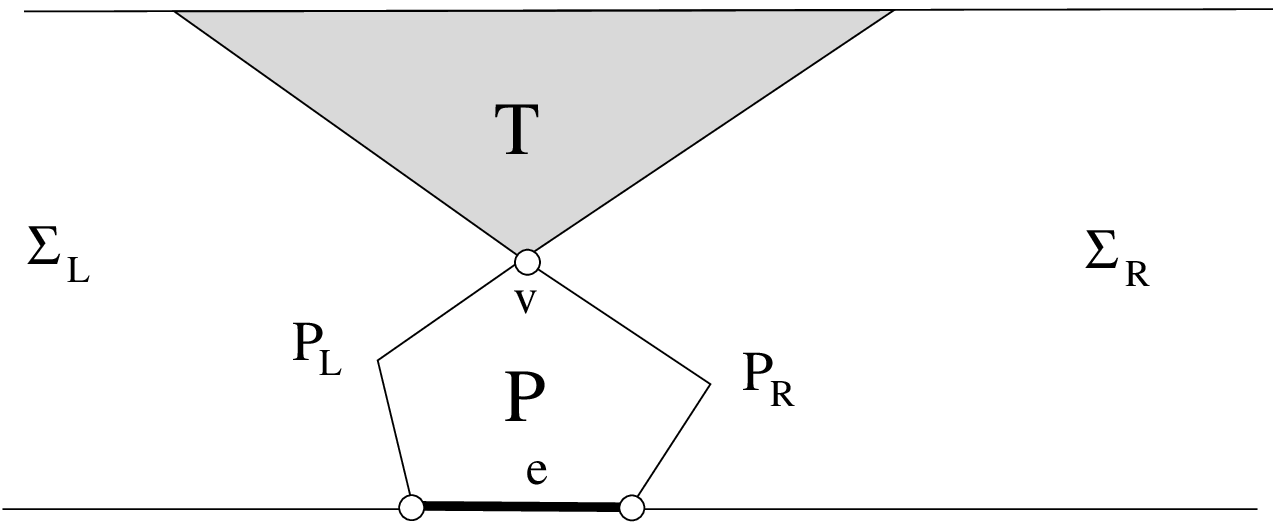}}
\newline
{\bf Figure 5.2:\/} Modified Strip
\end{center}
Let $P_L$ and $P_R$ repectively be the closures
of the left and right halves of $P-e-\{v\}$. Let
$\Sigma_L$ and $\Sigma_R$ denote the left and right
halves of $\Sigma-P-T$.

Let $\{p_i\}$ be an outer billiards orbit, with
$p_{i+2}=\psi(p_i)$.
 Let $q_i \in P$ be the midpoint
of the segment joining $p_i$ to $p_{i+1}$.  

\begin{lemma}
\label{tech1}
Suppose that $q_1,q_2 \in P_L$ and the clockwise arc from $q_2$ to $q_3$ does not
contain $e$.
Then $p_2 \in \Sigma_L$.
\end{lemma}

\startproof
Let $R$ be the
ray that starts at $v$ and moves up and to the left,
extending the left edge of $T$.
The conditions $q_1,q_2 \in P_L$ force $p_2$ to lie in the
region $\Sigma_L'$ bounded by $P_L$, by $R$, and by the negative $x$-axis.
$\Sigma_L$ is exactly the portion of $\Sigma_L'$ that lies below
the line $y=1$.  If $p_2 \in \Sigma_L'-\Sigma_L$ then $q_2=v$
and $p_3$ lies below the line $y=0$.  
 But then the clockwise arc from $q_2$ to $q_3$ contains
$e$.  This contradiction shows that $p_2 \in \Sigma_L$.
\endproof

\begin{lemma} 
\label{tech2}
Suppose that $q_1,q_2 \in P_R$ and the counterclockwise
arc from $q_1$ to $q_0$ does not contain $e$.
Then $p_2 \in \Sigma_R$. 
\end{lemma}

\startproof
This follows from the previous result and reflection symmetry.
\endproof

\subsection{The End of the Proof}

We label the vertices of $P$ so that $b<c$.   Let
\begin{equation}
p \in T_-(b) \cap T_+(c); \hskip 30 pt
j \in \{a,...,(b-1)\}.
\end{equation}
Our goal is to show that $p \in \Sigma_j$. 
 The sets of interest so us, such as $P_L$ and
$P_R$, all depend on an edge $e_j$, and we sometimes
write $P_L(e_j)$ and $P_R(e_j)$ to denote this
dependence.  We will establish the stronger result that

\begin{equation}
p \in \Sigma_L(e_j) \cup \Sigma_R(e_j).
\end{equation}

We set $p_2=p$ and we let $\{p_k\}$ be the
outer billiards orbit of $p$.  We let $q_k$ be the midpoint
of the line segment joining $p_k$ to $p_{k+1}$, as
in the previous section.
The next lemma refers to a definition made in \S \ref{harmony}.

\begin{lemma}
\label{harmony2}
The triple $(p_2,S_b,S_c)$ is a harmonious triple.
\end{lemma}

\startproof
The point $p_2$ lies in some tile $T_-(b,a)$.  This means
that $b$ is the starting point of an admissible path whose 
first spoke is $S_b$.  Since $(b,a)$ labels the tile
$T_-(b,a)$ in the backwards partition, the line
connecting $p_1$ to $p_2$ is tangent to $P$ at an
endpojnt of $S_b$.  In short $q_1$ is an endpoint
of $S_a$.  Similarly, $p_2$ is in some tile
$T_+(c,d)$.  The same argument shows that $q_2$
is an endpoint of $S_b$.
\endproof

We will treat the case when the spokes
$S_b$ and $S_c$ do not share a vertex,
though this is mainly for the sake of
drawing one picture rather than two.
When $S_b$ and $S_c$ share
a vertex, the proof is essentially
identical.

Let $X=S_b \cap S_c$.  In the case we
are considering, $X$ is a
point interior to both $S_b$ and $S_c$.
We rotate so that $S_b$ has positive
slope and $S_c$ has negative slope.

Let $\Sigma=\Sigma_j$ for some relevant index $j$.
Let $e_j$ be the edge of $P$ contained in
$\partial \Sigma_j$.  We can further rotate $P$,
if necessary, so that $e=e_j$ lies below the
intersection $S_b \cap S_c$.  Figure 5.3 shows
the situation.

Referring
to Lemma \ref{between}, we let $V(b,c,1)$ denote
the vertex set on the bottom of the picture
and $V(b,c,2)$ the vertex set on the top.
We will treat the case when $p_2$ lies to
the left of $X$.  This case relies on
Lemma \ref{tech1}.  The case when $p_2$
lies to the right of $X$ has practically
the same proof, and relies on Lemma \ref{tech2}
instead.

Let $H_b$ and $H_c$
respectively denote the right halfplanes
bounded by $S_b$ and $S_c$.  Figure 5.3 shows a
circle in place of a nice polygon.  This
is our attempt to draw a somewhat schematic
picture of what is going on. 

\begin{enumerate}
\item  Combining Lemma
\ref{harmony2} and Lemma \ref{subtend}, we
see that $q_1=S_b(1)$ and $q_2=S_c(2)$ both lie to the
left of $X$.
\item By Lemma \ref{between}, we have
$e_j \subset {\rm closure\/}(H_b-H_c)$.
\item By Lemma \ref{between} again, $v_j \subset H_c$. Here
$v_j$ is the vertex farthest from $e_j$.
\item By Steps 3 and 4, we have
$A(p_2) \subset P_L(e_j)$.  See \S \ref{harmony}
for a defintion of $A(P_2)$.
Hence $q_1,q_2 \in P_L(e_j)$.
\item The ordered pair $(q_2,q_3)$ is an admissible pair.
Hence, $q_3 \subset H_c$. 
\item Combining Steps 2 and 5, we conclude that
the clockwise arc connecting $q_2$ to $q_3$
does not contain $e_j$.
\item Combining Steps 4 and 6 and Lemma \ref{tech1},
we see that $p_2 \subset \Sigma_L(e_j)$, as desired.
This completes the proof.
\end{enumerate}

\begin{center}
\resizebox{!}{3.3in}{\includegraphics{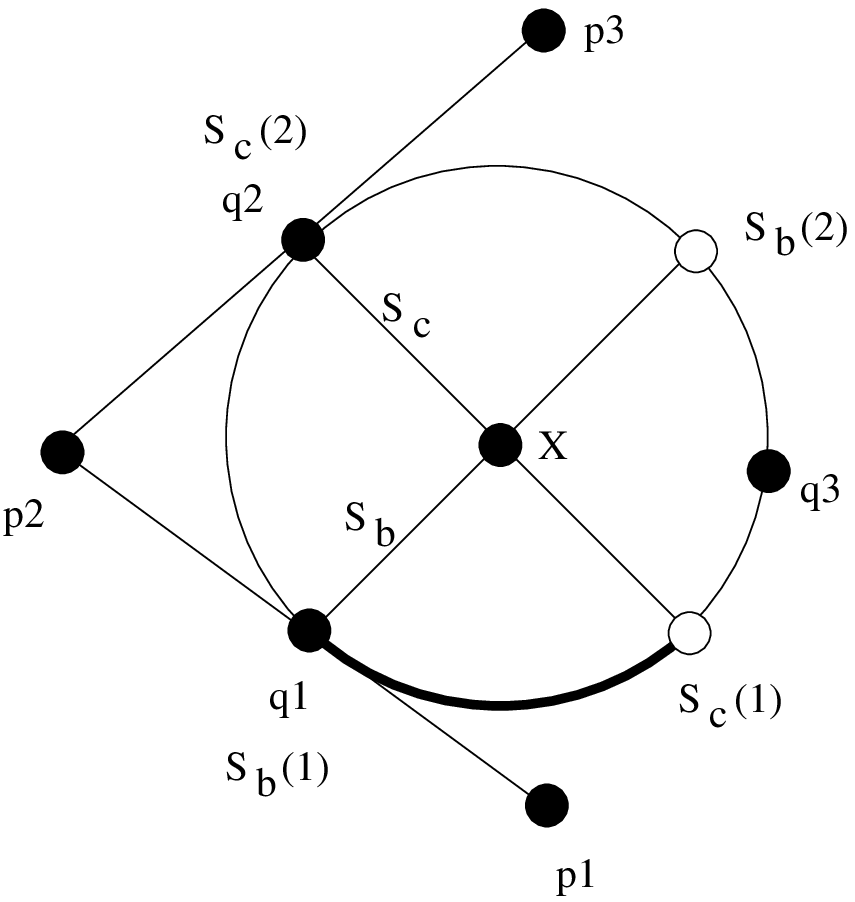}}
\newline
{\bf Figure 5.3:\/} Case 1
\end{center}

\newpage

\section{Theorem \ref{structure2} modulo a detail}
\label{proof2}

\subsection{Reformulation of the Result}

The unbounded regions correspond to the admissible paths
of length $1$.  When $a=b$, the result is a tautology.
So, we assume that $a \not = b$.  The admissible
path in question has length at least $3$.
We are going to prove Theorem \ref{structure2}
by induction.  
As previously, we replace $b$ by $b+n$ if necessary so
that $a<b$.

The admissible path $a \to b$ does not necessarily
involve all the spokes between $a$ and $b$.  
If $S_i$ is not involved in the path
$a \to b$, let $W_i=0$.  Otherwise,
let $W(a,b,i)$ be the vector that points
from the first endpoint of $S_i$ to the
last one.  By Lemma \ref{orientation},
we have $W_i=V_i$ for $i<b$.  For $i=b$
we have $W_i= \pm V_i$, where the
$(+)$ option is taken if and only if
$S_b$ is an ordinary spoke.  The
vector $W_i$ depends (mildly) on
$a \to b$, but we suppress this from our
notation.

\begin{lemma}
\label{move}
For any $p \in T(a \to b)$ we have
$$\psi(p)-p=\sum_{i=1}^b 2W_i.$$
\end{lemma}

\startproof
Let $(v,w)$ be the admissible pair of vertices
associated to $T(a \to b)$.  Recall that
$a$ and $b$ respectively name the first and
last spoke of the admissible path associated to
$T(a \to b)$ whereas $v$ and $w$ respectively
name the first endpoint of the path and the
last endpoint of the path.  The path $a \to b$
simply traces out the involved spokes.  By
definition
$$w-v=\sum_{i=a}^b W_i.$$
At the same time
$\psi(p)-p=2(w-v).$
Putting these two equations together gives the lemma.
\endproof

Lemma \ref{move} establishes an identity between
certain multiples of the vectors involved in the relevant
strip maps.  This is a good start.  What connects the
result in Lemma \ref{move} to the pinwheel map is the
claim that the multiples involved are precisely the
ones that arise in the relevant strip maps.  This
amounts to showing that certain translates of the
tile $T(a \to b)$ lie in the right place with respect to
the relevant strips.  Here is the main construction
in our proof.  Define

\begin{equation}
\label{tileversion}
T(a \to b;k)=T(a \to b)+\sum_{i=a}^k 2W_i; \hskip 30 pt
k=a,...,b
\end{equation}
The sets $T(a \to b;k)$ for $i=a,...,b$ are
translates of $T(a \to b)$.

Recall that $\mu_b$ is the local strip map with index $b$.
In this chapter we will prove
\begin{lemma}
\label{pin2}
$\mu_b(p)=p+2W_b$ for all $p \in T(a \to b,b-1)$.
\end{lemma}
In the next chapter, we will prove
\begin{lemma}
\label{pin1}
$T(a \to b;k) \subset \Sigma_k$ for all $k=a,...,(b-1)$.
\end{lemma}

\noindent
{\bf Remarks:\/} \newline
(i)  Lemma \ref{pin1} is the much
more interesting of the two results. 
It turns out that Lemma \ref{pin2}
is just a disguised version of
Corollary \ref{prepin3}. \newline
(ii)
It would be easier to say simply that
$T(a \to b;k) \subset \Sigma_k$ for all
$k=a,...,b$, but this is not true.
Lemma \ref{pin2} makes the strongest
statement we can make. 
\newline
\newline
\noindent
{\bf Proof of Theorem \ref{structure2}:\/}
We relabel so that $a=1$.
Let 
\begin{equation}
T_0=T(1 \to b); \hskip 40 pt
T_k=T(1 \to b;k).
\end{equation}
Let $p=p_0 \in T_0$ be an arbitrary point.
Define
\begin{equation}
\label{pointversion}
p_k=p_0+\sum_{i=1}^k 2W_i; \hskip 30 pt
k=1,...,b.
\end{equation}
This is just a pointwise version of
Equation \ref{tileversion}.

By Lemma \ref{move}, we have
\begin{equation}
\label{easy}
\psi(p_0)=p_b.
\end{equation}

Choose any $k=0,1,...,b-2$ and
consider the pair $(p_k,k)$. There are
two cases to consider.
Suppose first that the index $k$ is involved in 
the path $1 \to b$.   Then
\begin{equation}
p_k+2W_k \in \Sigma_k,
\end{equation}
by Lemma \ref{pin1}.
Therefore
\begin{equation}
\mu_{k+1}(p_k)=p_k+2W_k=p_{k+1}; \hskip 30 pt
\mu_{k+1}(p_{k+1})=p_{k+1}.
\end{equation}
Hence
\begin{equation}
(\psi^*)^2(p_k,k)=\psi^* \circ \mu_{k+1}(p_k,k)=
\psi^*(p_{k+1},k)=(p_{k+1},k+1).
\end{equation}
Now suppose that the index $k$ is not involved
in the path $1 \to b$.  Then $p_{k+1}=p_k$ and
this common point lies in $\Sigma_{k+1}$.
Hence
\begin{equation}
\psi^*(p_k,k)=(p_k,k+1)=(p_{k+1},k+1).
\end{equation}
In either case, we see that
\begin{equation}
(p_{k+1},k+1)=(\psi^*)^e(p_k,k)
\end{equation}
for some exponent $e=e_k$.
Applying this argument for as long as we can, we see that
\begin{equation}
(p_{b-1},b-1) = (\psi^*)^e(p_0,0),
\end{equation}
for some exponent $e$.
Finally, by  Lemma \ref{pin2},
we have
\begin{equation}
\psi^*(p_{b-1},b-1)=(\mu_b^*(p_{b-1}),b-1)=
(p_{b-1}+W_b,b-1)=(p_b,b-1).
\end{equation}
Hence $(p_b,b-1)$ is in the forward $\psi^*$-orbit
of $(p_0,0)$.  Combining this information with
Equation \ref{easy}, we see that
there is some positive $k<2n$ such that
\begin{equation}
(\psi(p),b-1)=(p_b,b-1)=(\psi^*)^k(p,0).
\end{equation}
This completes the proof.
\endproof

The rest of this chapter is devoted to the proof of
Lemma \ref{pin2}.  The material in the next section
will also be used in \S 7.

\subsection{The Conjugate Polygon}
\label{conjugate}

As in the proof of Theorem \ref{structure3}, we find it
useful here to consider both the forward and backward partition
of $P$.  Since our labelling conventions have been geared towards
the forward partition, we find it useful to set things up in a
way that doesn't require us to deal directly with the backward
partition. Let $R$ denote reflection in the $x$-axis.
For any subset $A \subset \R^2$, let
$\overline A=R(A)$.  One could equally well think of $R$ as
complex conjugation.

A basic and easy fact is that
$R$ maps the backward partition of $P$ to the
forward partition of $\overline P$.  For that reason,
we will consider the forward partitions of $P$ and
$\overline P$ at the same time. 
We rotate so that one edge of $P$ is horizontal and lies in the
$x$-axis.  We call this edge $e_0$.  The associated pinwheel
strip is $\Sigma_0$ and the associated spoke is $S_0$.
We also arrange that $P$ lies in the upper half plane,
so that $e_0$ is the bottom edge.  The labellings of
the remaining objects of $P$ are then forced by our
earlier conventions.  To emphasize the dependence
on $P$, we write $e_0(P)$, etc.
We let $e_0(\overline P)$ be the horizontal edge of $\overline P$.

The cyclic ordering forces two sets of equations.
\begin{equation}
\label{noshift}
\overline{e_k(P)} = e_{-k}(\overline P); \hskip 40 pt
\overline{\Sigma_k(P)} = \Sigma_{-k}(\overline P); \hskip 40 pt
\end{equation}

\begin{equation}
\label{shift}
\overline{S_k(P)}=S_{1-k}(\overline P); \hskip 40 pt
\overline{V_k(P)}=V_{1-k}(\overline P).
\end{equation}
This second set of equations is more subtle.
Figure 6.3 illustrates why
$\overline{S_0(P)}=S_{1}(\overline P).$
There are several possible pictures, depending on the
geometry of $P$, and we have picked one of the possiblities.
The other possibilities have the same outcome.
The rest of Equation \ref{shift} is then forced by
the cyclic ordering.

\begin{center}
\resizebox{!}{1.9in}{\includegraphics{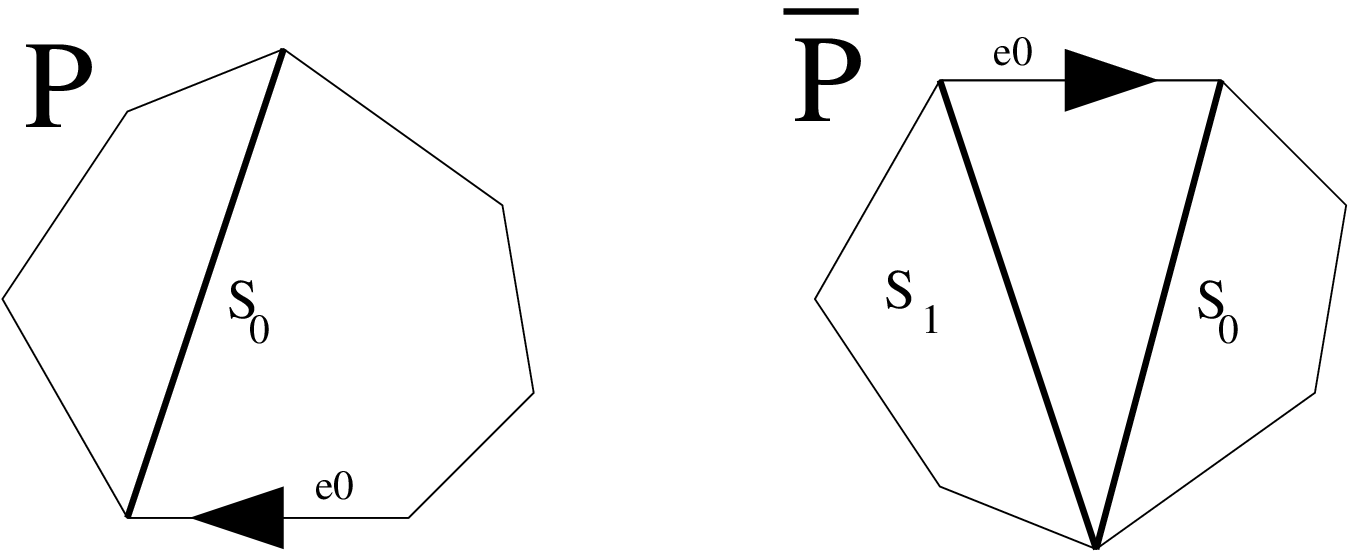}}
\newline
{\bf Figure 6.3:\/} The polygon and its conjugate
\end{center}

We can define the sets $T(a \to b;k)$ relative to $\overline P$
just as well as for $P$.  We tack a $P$ or $\overline P$ on the
end of our notation to indicate which polygon we mean.

\begin{lemma}
\label{carry}
$\overline{T(a \to b; k;P)}=T(1-b,1-a;-k;\overline P).$
\end{lemma}

\startproof
We've already remarked that $\psi$ sets up a bijection between
the tiles in the forward partition of $P$ to the tiles in the
backward partition.  Briefly using the notation from the
previous chapter, we have
\begin{equation}
\label{oldnotation}
\psi(T_+(a \to b))=T_-(b \to a).
\end{equation}
The path $b \to a$ is the same as the path $a \to b$
but it is given the opposite orientation.  We call
this the {\it reversal property\/}.  We will use it
below.

The composition $R \circ \psi$ carries the tiles in the
forward partition of $P$ to the tiles in the forward
partition of $\overline P$.  
Combining Equations \ref{shift} and \ref{oldnotation}, we see that
\begin{equation}
\label{firstreflect}
R \circ \psi(T(a \to b;P)) = T(1-b \to 1-a;\overline P).
\end{equation}

Note that
\begin{equation}
\label{obvious}
\psi(T(a \to b;P))=T(a \to b;b;P).
\end{equation}
Combining Equations \ref{obvious} and \ref{firstreflect}, we have
\begin{equation}
\overline{T(a \to b;b;P)}= T(1-b,1-a;\overline P).
\end{equation}
By the reversal property and
Equation \ref{shift}, we have
\begin{equation}
\overline{W_k(P)}=-W_{1-k}(\overline P).
\end{equation}

To consider the case $k=b$ in detail, we have
$$\overline{T(a \to b;b-1;P)} =
\overline{T(a \to b;b;P)} - \overline{W_b(P)}=$$
$$T(1-b \to 1-a;\overline P)+W_{1-b}(\overline P)=
T(1-b,1-a;1-b;\overline P).$$
in short
\begin{equation}
\label{prepin4}
\overline{T(a \to b;b-1;P)}=T(1-b,1-a;1-b;\overline P).
\end{equation}
Applying the same argument, inductively, to each of the
vectors $W_{b-1}(P)$, $W_{b-2}(P)$, etc., we establish the lemma.
\endproof

\subsection{Proof of Lemma \ref{pin2}}

We are going to apply Corollary \ref{prepin3}
to $\overline P$.
There are two cases for us to consider, depending on whether
or not the spoke $S_b$ is ordinary.  We'll first consider
the case when $S_b$ is ordinary.

It is convenient to set 
\begin{equation}
\alpha=1-b; \hskip 40 pt \beta=1-a.
\end{equation}
By Lemma \ref{carry},
\begin{equation}
\overline{T(a \to b,b-1;P)}=T(\alpha \to \beta;\alpha;\overline P).
\end{equation}
By definition and Lemma \ref{orientation},
\begin{equation}
T(\alpha \to \beta;\alpha;\overline P)=
T(\alpha \to \beta;\overline P) + V_{\alpha}(\overline P).
\end{equation}
By Corollary \ref{prepin3}, the line
$$\partial_1 \Sigma_{\alpha-1}(\overline P)$$
separates $T(\alpha \to \beta;\overline P)+V_{\alpha}(\overline P)$
from $\overline P$.  Applying the reflection $R$ and using Equation \ref{shift},
we see that 
the line
$$\partial_1 \Sigma_{b}(P)$$ 
separates $T(a,b;b-1;P)$ from $P$.
But then $$\mu_b(p)=p+V_b$$ for any $p \in T(1 \to b;b-1;P)$.
Since $S_b$ is an ordinary spoke, $V_b=W_b$ by Lemma \ref{orientation}.
Now we know that $\mu_b(p)=p+W_b$, as desired.

When $S_b$ is special, the proof is the same except for some sign
changes.  This time, the line
$$\partial_2 \Sigma_{b}(P)$$ 
separates $T(a,b;b-1;P)$ from $P$.
But then $$\mu_b(p)=p-V_b$$ for any $p \in T(1 \to b;b-1;P)$.
This time $-V_b=W_b$, and we get the same result as in the
previous case.

This completes the proof of Lemma \ref{pin2}.

\newpage

\section{Proof of Lemma \ref{pin1}}

\subsection{Some Combinatorial Definitions}

\noindent
{\bf Abstract Admissible Paths:\/}
We say that an {\it abstract admissible path\/} is a finite
tree $\tau$ with the following structure.  First, there is a
distinguished maximal path $\gamma$ in $\tau$ having odd length
at least $3$.
Every other edge of $\tau$ is incident to $\gamma$.
We call the edges of $\tau-\gamma$ {\it special\/}.
We call the edges of $\gamma$ {\it ordinary\/},
except perhaps for the first and last edge
of $\gamma$.  The first and last edges of
$\gamma$ can be either special or ordinary.
We draw $\gamma$ as a zig-zag with lines
of alternating negative and positive
slope.  We insist that every edge of
$\tau$ intersects $x$-axis.

We orient $\gamma$ from left to
right.  We draw the ordinary edges with
thick lines and the special edges with 
thin lines. Figure 7.1 shows an example.
A dotted line represents the $x$-axis.
$\gamma$ is the path $13678$.
Here, the first edge of
$\gamma$ is ordinary and the last
edge is special. The {\it initial
vertex\/} is the left endpoint of the
$\gamma$.  The {\it final edge\/} if
the last edge of $\gamma$.

\begin{center}
\resizebox{!}{1.7in}{\includegraphics{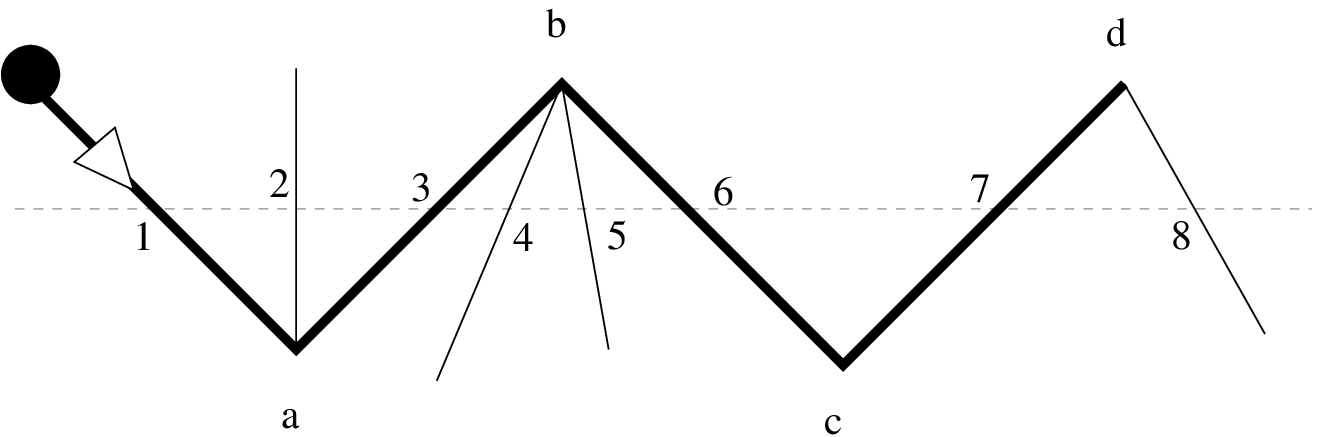}}
\newline
{\bf Figure 7.1:\/} Abstract admissible paths.
\end{center}

\noindent
{\bf Linear Order:\/}
There is a natural linear ordering on
the edges of $\tau$, induced from the order
in which they intersect the $x$-axis,
goind from left to right.  The
numerical labels in Figure 7.1 indicate the ordering.
We see that $\tau' \subset \tau$ is a {\it prefix\/}
of $\tau$ of $\tau'$ and $\tau$ share the same initial
set of edges and if $\tau'$ is an abstract
admissible path in its own right.
\newline
\newline
{\bf Flags and Sites:\/}
To each edge $e$ of $\tau$, except the last one,
we assign a vertex $v_e$.  If $e$ is an edge
of $\gamma$, then $v_e$ is the leading vertex
of $e$.  If $e$ is not an edge of $\gamma$, then
$v_e$ is the vertex of $\gamma$ incident to
$e$.  We say that a {\it flag\/} is a pair
$(e,v_e)$.  For instance, the flags in
Figure 7.1 are
$$(1,a); \hskip 5 pt
  (2,a); \hskip 5 pt
  (3,b); \hskip 5 pt
  (4,b); \hskip 5 pt
  (5,b); \hskip 5 pt
  (6,c); \hskip 5 pt
  (7,d).
$$
We say that a {\it site\/} is a pair $(\tau,f)$, where
$f$ is a flag of $\tau$.
\newline
\newline
{\bf Natural Involution:\/}
There is a natural involution $R$ on the set of abstract
admissible paths:  Simply rotate the path about the
origin by $180$ degrees and you get another one.
We call this map $R$.   The map $R$ carries flags
of $\tau$ to flags of $R(\tau)$ in a slightly
nontrivial way.  We first create {\it reverse flags\/}
of $\tau$ by interchanging the notion of left and
right, and then we apply $R$ to these reverse flags
to get ordinary flags of $R(\tau)$.
In Figure 7.1, the reverse flags are
$$(8,d); \hskip 5 pt
  (7,c); \hskip 5 pt
  (6,b); \hskip 5 pt
  (5,b); \hskip 5 pt
  (4,b); \hskip 5 pt
  (3,a); \hskip 5 pt
  (2,a).
$$
$R$ maps the leftmost flag of $\tau$ to the
image under rotation of the leftmost
reverse flag. For instance $(1,a)$ corresponds
to the rotation of $(2,a)$.
\newline
\newline
{\bf Reduction:\/}
We say that the site $(\tau',f)$ is a {\it direct reduction\/}
of $(\tau,f)$ if $\tau'$ is a prefix of $\tau$.  The flag
$f$ is the same in both cases.  We say that
$(\tau',f')$ is an {\it indirect reduction\/} of
$(\tau,f)$ if $(\tau',f')$ is a direct reduction
of $R(\tau,f)$.  We say that one site $(\tau_2,f_2)$
a {\it reduction\/} of another site $(\tau_1,f_1)$ if
$(\tau_2,f_2)$ is either a direct or an indirect
reduction of $(\tau_1,f_1)$.  In this case, we write
$(\tau_1,f_1) \to (\tau_2,f_2)$.
\newline
\newline
{\bf Hereditary Properties:\/}

Let ${\cal C\/}$ be a collection of sites.  We say
that ${\cal C\/}$ is {\it hereditary\/} if
it has the following properties.
\begin{itemize}
\item $\cal C$ is closed under the natural involution.
\item $\cal C$ is closed under reduction.
\end{itemize}

Say that a site $(\tau,f)$ is
{\it prototypical\/} if $f$ is the first flag of $\tau$.
 Let
$\Omega$ be a map from $\{0,1\}$.  We say that $\Omega$ is
{\it hereditary\/} if $\Omega$ has
the following properties.
\begin{itemize}
\item $\Omega$ evaluates to $1$ on all prototypical sites in $\cal C$.
\item $\Omega \circ R=\Omega$. Here $R$ is the natural involution.
\item If $\Omega(\tau_1,f_1)=1$ and
$(\tau_2,f_2) \to (\tau_1,f_1)$ then
$\Omega(\tau_2,f_2)=1$.
\end{itemize}

\subsection{The Reduction Lemma}

In this section we will prove the following result.
\begin{lemma}[Reduction]
Suppose that $\cal C$ is a hereditary collection
of sites and $\Omega$ is a heriditary function on
$\cal C$.  Then $\Omega \equiv 1$ on $\cal C$.
\end{lemma}

\startproof
It suffices to prove that, through the
two operations of $R$ and reduction, every
site can be transformed into a prototypical site.
It henceforth 
goes without saying that all sites belong
to $\cal C$.

Let $(\tau,f)$ be a site.  Let $\gamma$ be
the maximal path of $\tau$.  Let $f=(e,v)$.
Either $v$ lies in the left half of $\gamma$
or the left half. (There are an even number
of vertices.)  Applying $R$ if necessary,
we can assume that $v$ lies in the left
half of $\tau$.  If $\gamma$ has length
$5$ we let $\tau'$ denote the subtree
of $\tau$ obtained by deleting the last
two vertices of $\gamma$ and all incident
edges.  Then $\tau'$ is a prefix of
$\tau$ and $(\tau,f) \to (\tau',f)$.

\begin{center}
\resizebox{!}{1.3in}{\includegraphics{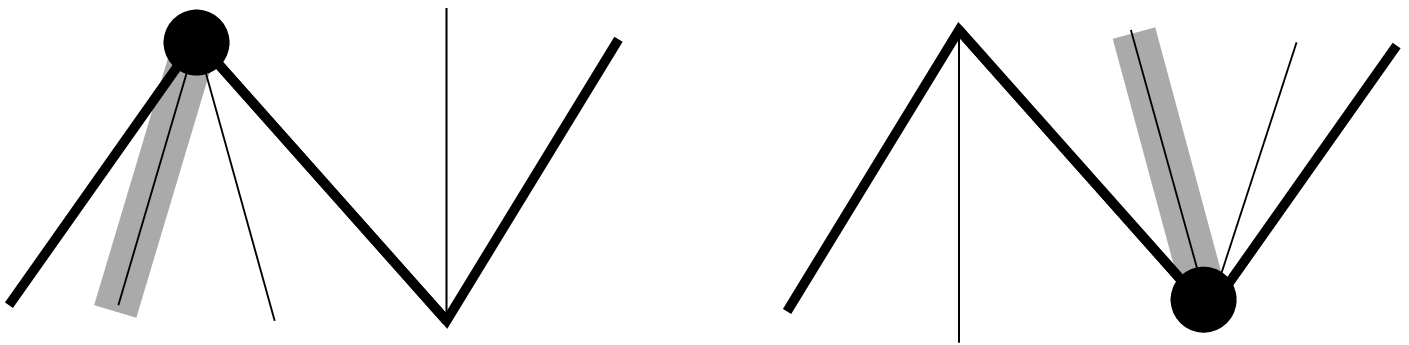}}
\newline
{\bf Figure 7.2:\/} Abstract admissible paths.
\end{center}

We just have to worry about the case when $\gamma$ has
length $3$.  Let $(\tau_1,f_1)=(\tau,f)$ and
let $(\tau_2,f_2)=R(\tau,f)$.  Also,
let $e_k$ be the edge of $f_k$ for $k=1,2$.
If $e_1$ is not the first edge of $\tau_1$
then $(\tau_2,f_2)$ has the
following two properties.
\begin{enumerate}
\item $e_2$ is neither of the last two edges of
$\tau_2$.
\item At least $3$ edges of $\tau_2$ are incident
to the right vertex of $\gamma$.
\end{enumerate}
Figure 7.2 shows a typical situation.
The thick grey lines represent $e_1$ and $e_2$.
The upshot is that after applying $R$, we can
assume that $f=(e,v)$, where $e$ is neither of the
last two edges of $\tau$.  We let $\tau'$ be
the prefix obtained by cutting off these last
two edges.  The second property mentioned above
guarantees that $\tau'$ is a prefix of $\tau$.
Again we have $(\tau,f) \to (\tau',f)$.

In summary, the process above only stops when
we reach a prototypical site.
\endproof

\subsection{The Pinwheel Collection}

Let $P$ be a nice polygon.  Each admissible path
associated to $P$ gives rise to an abstract
admissible path.  
This path encodes the way the spokes in the
path (and the skipped spokes) meet at their
endpoints.   Figure 7.3 shows an example.
The unusual labelling on the abstract
admissible path is present only to
underscore the correspondence.

\begin{center}
\resizebox{!}{2.4in}{\includegraphics{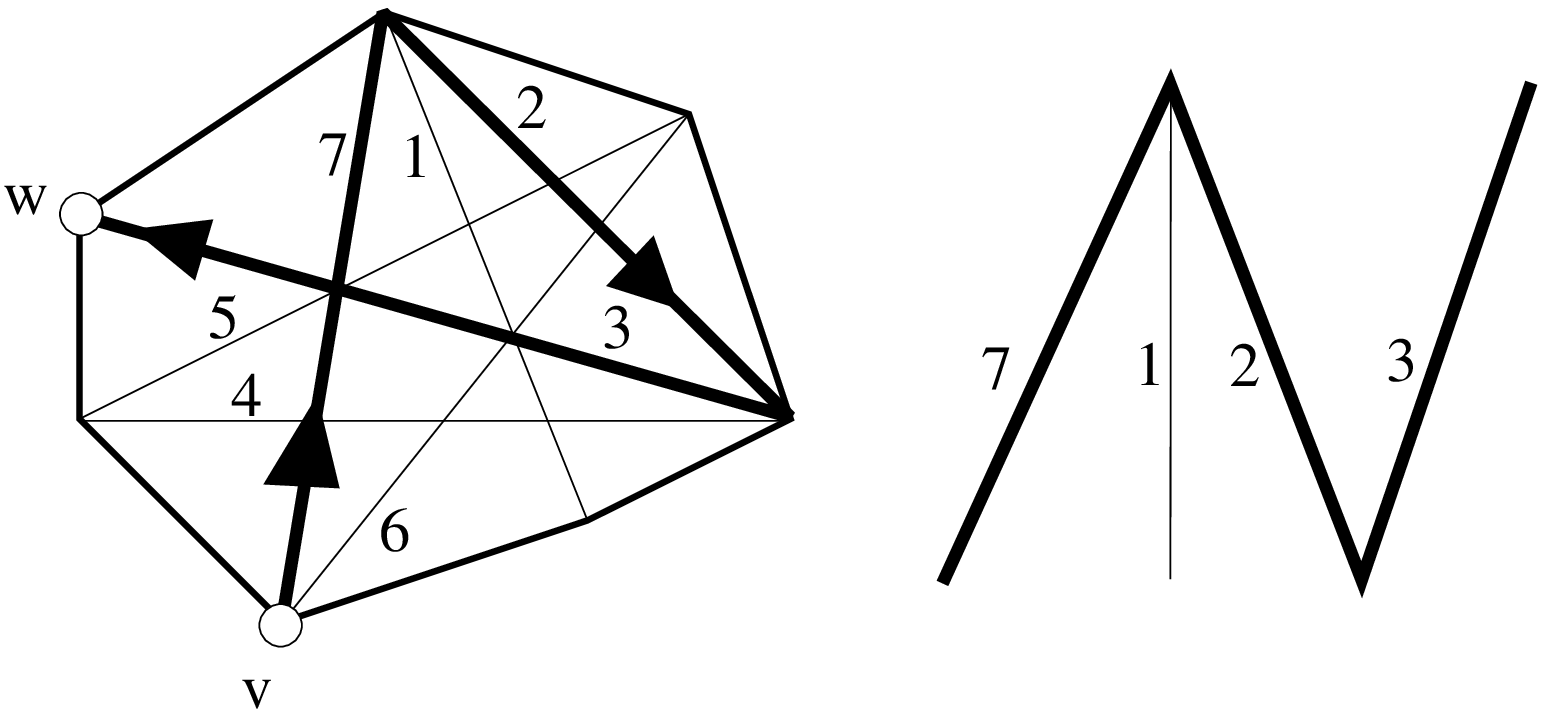}}
\newline
{\bf Figure 7.3:\/} An admissible path and its abstraction
\end{center}

We call the abstract admissible path
produced in this way the {\it abstraction\/}
of the admissible path.  The maximal
path in the abstract admissible path
corresponds to the actual admissible path.
We let $\cal C$ be the class of all sites
$(\tau,f)$, where $\tau$ is the abstraction of
an admissible path for some nice polygon.

\begin{lemma}
$\cal C$ is a heriditary collection.
\end{lemma}

\startproof
The fact that $\cal C$ is closed under reduction 
comes from the fact that we have built in the
basic properties of admissible paths into
our definition of abstract admissible paths.
If $\tau$ is the abstraction of $a \to b$,
then any prefix $\tau'$ is the abstraction
of some admissible path $a \to b'$, where
$b'<b$.

The analysis in \S \ref{conjugate} shows that
$\cal C$ is closed under the natural involution.
The basic reason, as we have discussed already,
is that $\psi$ carries the forward partition
to the backward partition. $\psi$ maps
the forward tile associated to the
path $a \to b$ to the backward tile
associated to the path $b \to a$.
The abstraction of $b \to a$ is exactly
the image of the abstraction of $a \to b$
under the natural involution.
\endproof

\subsection{The Binary Function}

Now we are going to define a function
$\Omega: {\cal C\/} \to \{0,1\}$.
We will first consider the situation for
a given nice polygon $P$, and then we
will take into account all nice
polygons at the same time.

Let $(\tau,f)$ be a site in $\cal C$.  This means
that there is a nice polygon $P$ and an
admissible path $a \to b$ such that
$\tau$ is the abstraction of $a \to b$.
Moreover, $f$ is just one of the sites.
The edges of $\tau$ are naturally in
correspondence with the strips
$\Sigma_a,...,\Sigma_b$.  Moreover, there
is a natural correspondence between
the sets
$T(a \to b;a),...,T(a \to b,b-1)$ and the
sites of $\tau$.  The correspondence
is set up in such a way that each
site $(\tau,f)$ corresponds to a pair
\begin{equation}
\label{pairs}
T(a \to b;k); \hskip 40 pt \Sigma_k.
\end{equation}
All of this depends on $P$.
We define
$\Omega'(\tau,f;P)=1$ if Lemma \ref{pin1} is
true for the pair in Equation \ref{pairs}
Finally, we define
$\Omega(\tau,f)=1$ if and only if
$\Omega'(\tau,f;P)=1$ for every instance
in which the site $(\tau,f)$ arises.

Lemma \ref{pin1} is equivalent to the statement that
$\Omega \equiv 1$ on $\cal C$. Accordingly, we will
prove Lemma \ref{pin1} by showing that $\Omega$
is a heriditary function and then invoking the
Hereditary Lemma.

\subsection{Property 1}

The prototypical sites correspond to the
case $k=a$ in Equation \ref{pairs}.
Cyclically relabelling, we take $a=1$.
The following lemma implies that $\Omega=1$
on all prototypical sites.

\begin{lemma}
\label{base case}
$T(1 \to b;1) \subset \Sigma_1$.
\end{lemma}

\startproof
We rotate so that $\Sigma_1$ is horizontal, and
$T(1 \to b)$ lies beneath the lower boundary
of $\Sigma_1$.  Let $L\Sigma_1$ and
$U\Sigma_1$ denote the lower and upper boundaries
of $\Sigma_1$ respectively.
$T(1 \to b)$ is contained 
in the shaded region $T$ shown in Figure 4.3.
By Lemma \ref{small}, every point of $T(1 \to b)$
is closer to $L\Sigma_1$ than (half) the distance
between $L\Sigma_1$ and $R\Sigma_1$.
We have the following $2$ properties.
\begin{itemize}
\item The vector $W_1$ points from $L\Sigma_1$
to $U\Sigma_1$.
\item $L\Sigma_1$ contains the top edge of $T(1 \to b)$.
\end{itemize}
It follows from these two properties that
$T(1 \to b)+W_1 \subset \Sigma_1$.
\endproof

\subsection{Property 2}

That $\Omega \circ R=\Omega$ is a consequence of the relations
between $P$ and $\overline P$ worked out in \S \ref{conjugate}.
Let $(\tau_1,f_1)$ be the site corresponding to the pair
in Equation \ref{pairs}.
We label the edges of $\tau_1$ as $a,...,b$.
Let $(\tau_2,f_2)=R(\tau_1,f_1)$.  We label the
edges of $\tau_2$ as $(1-b),...,(1-a)$.

We label the sites of $\tau_1$ by symbols of the
form $\langle k \rangle_1$.  Here $k$ names the
label of the edge involved in the site.  Likewise
we label the sites of $\tau_2$ with the label
$\langle k \rangle_2$. 

\begin{lemma}
 The natural involution
carries $\langle k \rangle_1$ to
$\langle -k \rangle_2$. 
\end{lemma}

\startproof
Given that the natural involution reverses the order of the sites,
it suffices to check the claim for a single site.  We choose
the first site, with $k=a$.  The natural involution carries
this site to the that involves the next-to-last edge of
$\tau_2$.  But this edge is labelled $-a$.
\endproof

According to this result, if the site $(\tau_1,f_1)$ corresponds
to the pair
$$
T(a \to b;k); \hskip 40 pt \Sigma_k.
$$
then the site $(\tau_2,f_2)$ corresponds to the pair
$$T(1-b \to 1-a;-k); \hskip 40 pt \Sigma_{-k}.$$
But, by Lemma \ref{carry}, reflection in the
$x$-axis carries the one pair to the other.
Hence, the desired containment holds in the one
case if and only it holds in the other.
In other words $\Omega(\tau_1,f_1)=\Omega(\tau_2,f_2)$.
This establishes the second property.

\subsection{Property 3}

This property has the most interesting proof.

Recall from \S \ref{tile geometry} that the primary
walls of the forward partition divide $\R^2-P$ into
primary cones.  The beginning vertex $v$ of the admissible
path $a \to b$ is the apex of the cone.
If we consider all the tiles of the form $T(a \to b)$
with $a$ fixed and $b$ increasing, then Lemma
\ref{pin1} involves increasingly many
containments.  On the other hand, the tiles
involved get smaller and smaller in the
sense that the shrink down to $v$.
One might expect that the vertex $v$ itself
satisfies all the identities we can form.

Let $a \to B$ denote the maximal admissible
path that starts with $S_a$.  Then $a \to B$
corresponds to the triangular tile $T(a \to B)$
that $v$ as a vertex.   
Define
\begin{equation}
p_0=v; \hskip 40 pt
p_k=p_0+\sum_{i=1}^k 2W_i.
\end{equation}
Note the similarity between this equation and
Equation \ref{tileversion}.
In the next section we will prove the
following result. Assume this result for now.	

\begin{lemma}
\label{apex}
$p_k \subset \Sigma_k$ for $k=1,...,B$.
\end{lemma}

When $k \leq b-1$ we
let $\Omega(a \to b;k)$ be statement that
$T(a \to b; k) \subset \Sigma_k$.

\begin{corollary}
\label{induct}
Suppose $c>b$ is that that $a \to c$ is an admissible
path.  Then $\Omega(a \to b;k)$ implies
$\Omega(a \to c,k)$.
\end{corollary}

\startproof
The admissible paths $a \to b$ and $a \to c$ agree 
except for possibly the last edge of $a \to b$.
Hence the sets
$$
T(a \to b;k); \hskip 30 pt
T(a \to c;k); \hskip 30 pt
p_k
$$
are respectively translates, by the same vector $X$, of the sets
$$
T(a \to b); \hskip 30 pt
T(a \to c); \hskip 30 pt
p_0.
$$
Let $\widehat T(a \to b)$ denote the convex hull
of $T(a \to b)$ and $p_0$.  Let 
$\widehat T(a \to b;k)$ denote the convex hull
of $T(a \to b;k)$ and $p_k$.  First of all,
\begin{equation}
T(a \to c) \subset \widehat T(a \to b),
\end{equation}
by the analysis done in connection with Figure 4.2.
Translating the whole picture by $X$, we have
\begin{equation}
T(a \to c;k) \subset \widehat T(a \to b;k) \subset \Sigma_k.
\end{equation}
The second equality follows from the convexity of
$\Sigma_k$.
\endproof

Corollary \ref{induct} is just a restatement of Property 3.

\subsection{Proof of Lemma \ref{apex}}

Now we take care of the final piece of business.

We relabel so that $a=1$.
Let $v_0=v=p_0$ and let $\{v_k\}$ be the successive
vertices of our admissible path, where we only
advance the point if the spoke is actually
involved in the path. Put another way,
the vertex $v_k$ is incident to the spokes
$S_k$ and $S_{k+1}$. Figure 7.4 shows an example.

\begin{center}
\resizebox{!}{2.7in}{\includegraphics{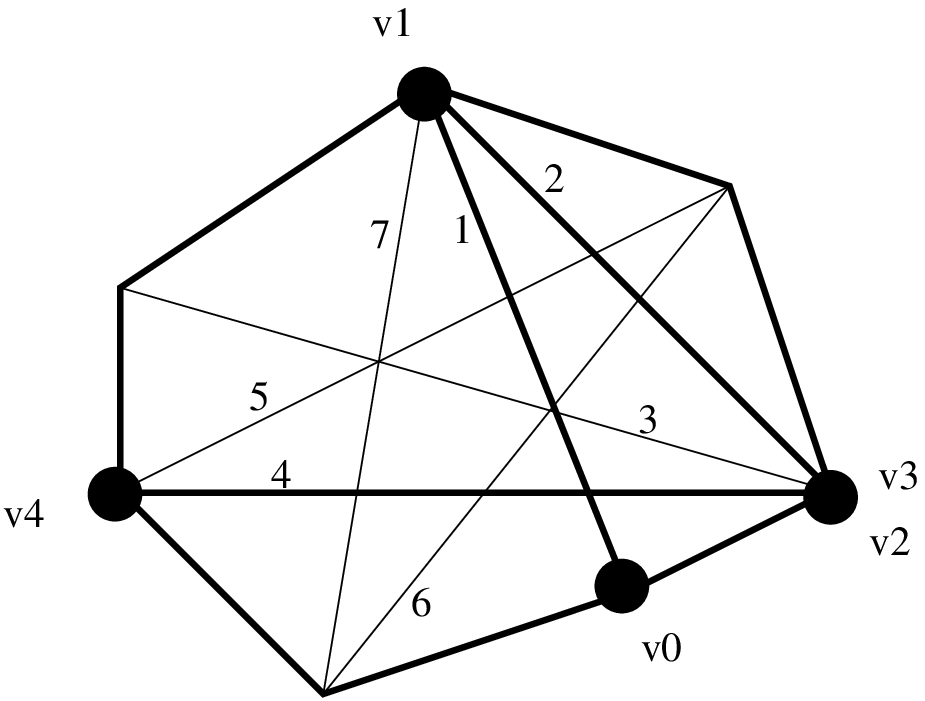}}
\newline
{\bf Figure 7.4:\/} The vertices of the path
\end{center}

The strip $\Sigma_k$ contains the three points of
$\partial S_k \cup \partial S_{k+1}$ and
$v_k=S_k \cap S_{k+1}$.
Compare Lemma \ref{minimal1}.
From this structure, we see that $v_k$ lies on the
centerline of the strip $\Sigma_k$ for all $k$.

\begin{center}
\resizebox{!}{1.4in}{\includegraphics{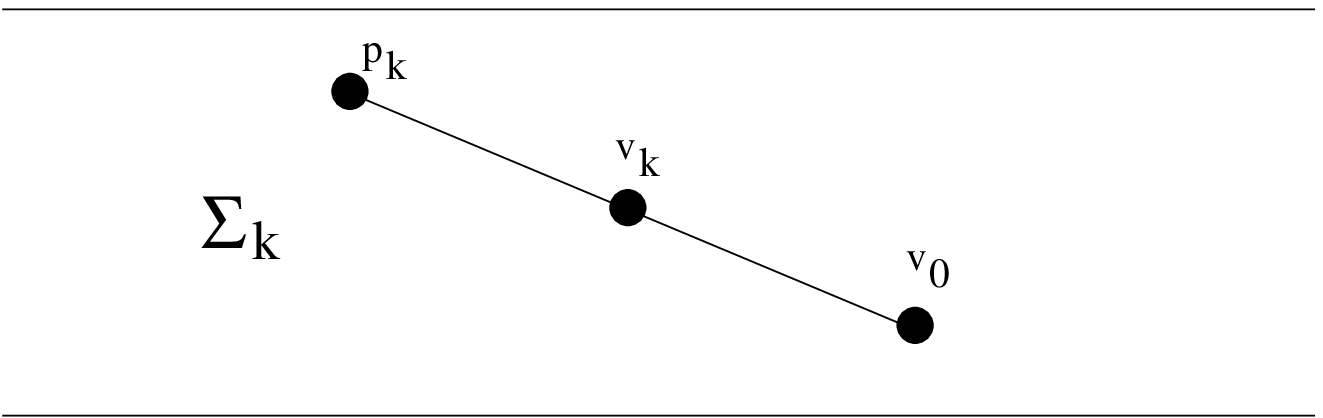}}
\newline
{\bf Figure 7.5:\/} The vertices of the path
\end{center}

The fixed point of $H$, namely $v_0$, belongs
to every strip.  In particular, $v_0 \in \Sigma_k$.
This situation forces $p_k=H(v_k) \in \Sigma_k$.
\endproof

\newpage

\section{Quasirational Polygons}
\label{qr}

\subsection{Notation}

Let $P$ be a quasirational $n$-gon.    We can scale
$P$ so that the quantities
\begin{equation}
A_j={\rm area\/}(\Sigma_j \cap \Sigma_{j+1})
\end{equation}
are all integers.  Let $D$ be the least common multiple
of $\{A_1,...,A_n\}$.  
Define
\begin{equation}
D_j=\frac{D}{A_j}.
\end{equation}
The quantities $D_1,...,D_n$ are all integers.

Let $e_j$ be the edge of $P$ contained in
the boundary of $\Sigma_j$.   As usual, we orient
$e_j$ so that, taken together, all these edges go
clockwise around $P$.
Let $W_j$ be the vector parallel to $e_j$ that
spans $\Sigma_{j+1}$.

The vector $V_{j+1}$ associated to the
strip $\Sigma_{j+1}$ points from one
corner of $\Sigma_j \cap \Sigma_{j+1}$ to the opposite
corner.      The point is that
the spoke associated to $\Sigma_{j+1}$  joins the
head point of $e_{j+1}$ to the tail point of $e_j$.   

We draw a picture in the case $j=1$.
The whole construction is affinely invariant,
so to draw pictures we normalize
so that $\Sigma_1$ is horizontal and
$\Sigma_2$ is vertical, and both have
the same width.
Up to dihedral symmetry, there 
are two possible local pictures.
These are shown in Figure 8.1.
Which case occurs depends on whether
$e_j$ and $e_{j+1}$ are adjacent.

\begin{center}
\resizebox{!}{2.5in}{\includegraphics{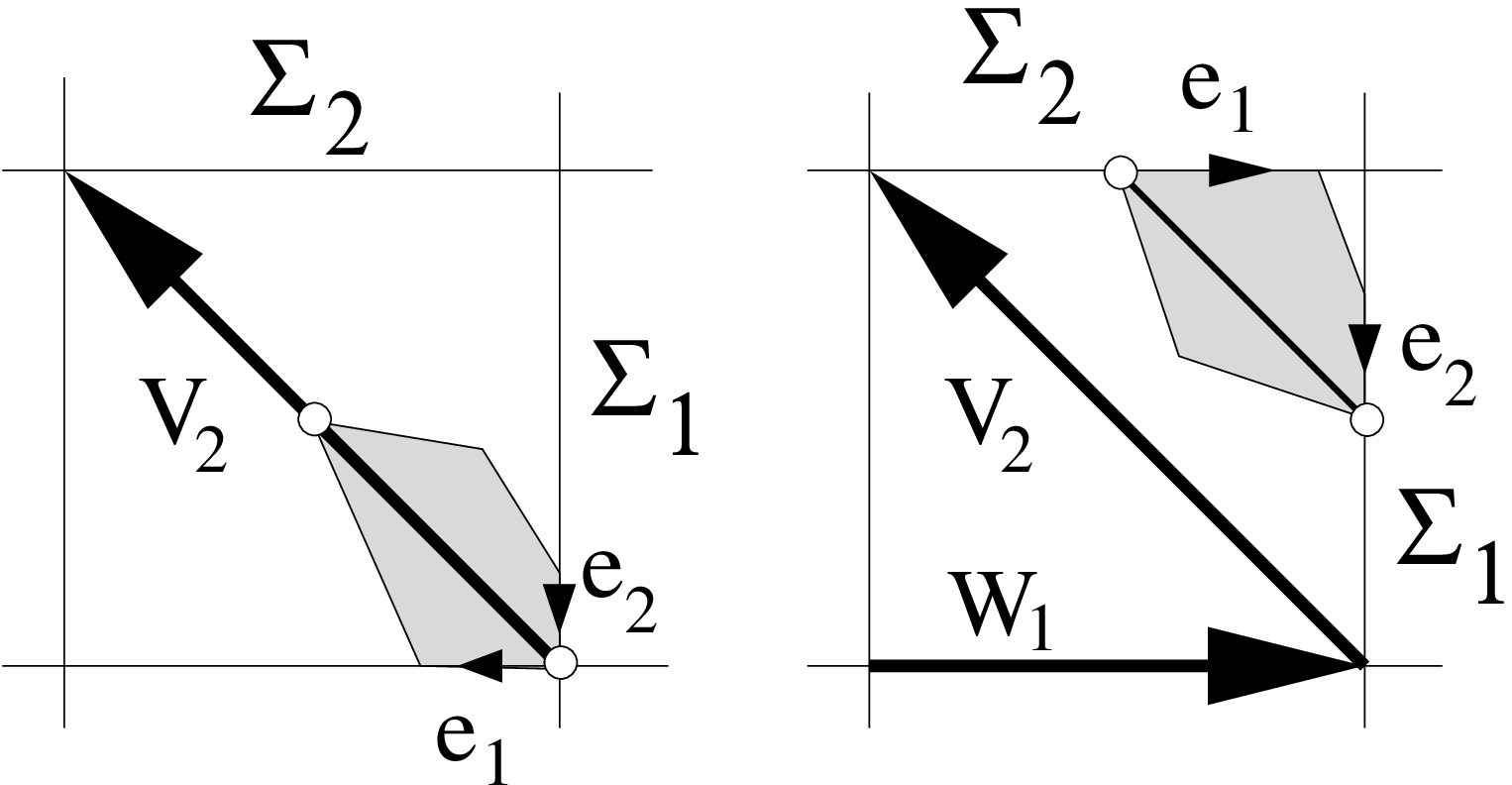}}
\newline
{\bf Figure 8.1:\/} Two consecutive strips
\end{center}

We will carry out
the construction in case $e_j$ and $e_{j+1}$ are
not adjacent. The other case has essentially the
same treatment.

As in Equation \ref{X}, define
\begin{equation}
\label{X2}
\widehat X=\bigcup_{j=1}^n (\Sigma_j \times \{j\}) \subset \R_n^2.
\end{equation}

Let $P_j=P \times \{j\}$.   The polygons
$P_1,...,P_n$ are just copies of $P$.      
Let $v_j$ be the vertex of $P_j$ that lies on the centerline
of $\Sigma_j \times \{j\}$.
Let $Q_j$ be the polygon obtained by rotating
$P_j$ by $180$ degrees about $v_j$.
Figure 8.2 shows a picture of these polygons
for $j=1,2$.  Technically, these polygons do not
live in the plane, but rather in $\R^2_n$.  However,
we have drawn their projections into the plane.
The polygons $P_j$ all have the same image in the
plane, so we set $P=P_1=P_2$ in the picture.

\begin{center}
\resizebox{!}{4in}{\includegraphics{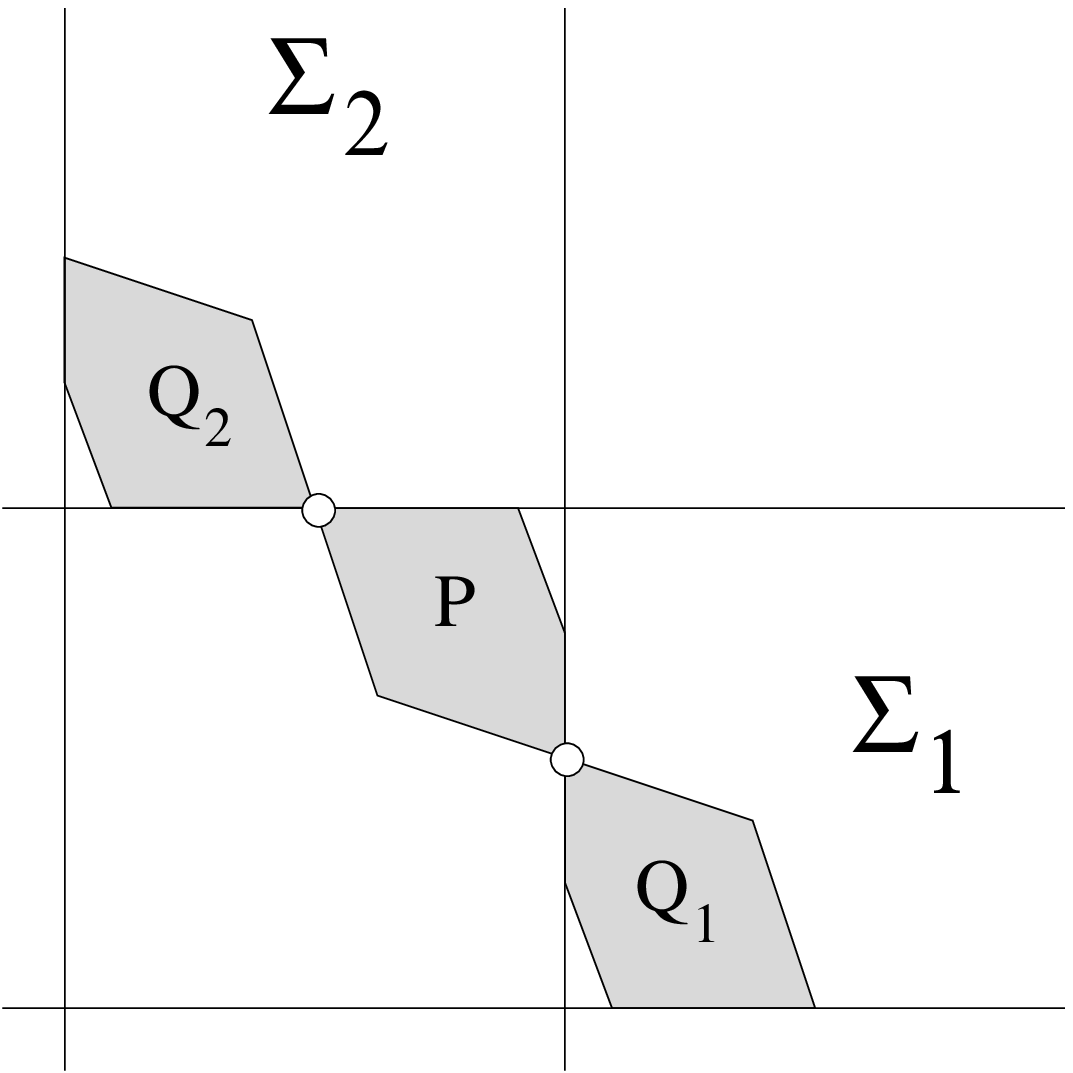}}
\newline
{\bf Figure 8.2:\/} $P$ and $Q_1$ and $Q_2$.
\end{center}

\subsection{The Main Argument}

We define
\begin{equation}
R_j^m = P_j^m \cup Q_j^m; \hskip 30 pt
P_j^m = P_j + (mW_j,j); \hskip 20 pt
Q_j^m = Q_j + (mW_j,j).
\end{equation}
We think of these as open polygons.

\begin{lemma}
\label{quasi1}
Let $m$ be any integer.   Then
$$\psi(R_j^{mD_j})=R_{j+1}^{mD_{j+1}}.$$ Moreover,
suppose $p \in \Sigma_j \times \{j\}$ lies between
$R_j$ and $R_j^{mD_j}$.   Then
$\psi(p)$ lies between $R_{j+1}$ and
$R_{j+1}^{mD_{j+1}}$.  
\end{lemma}

\startproof
We first consider the strips $\Sigma_1$ and $\Sigma_2$ in detail.
We normalize by an affine transformation so that the picture 
looks as in Figures 8.1 and 8.2.  Our two strips (as normalized)
are tiled by squares -- i.e., translates of
$\Sigma_1 \cap \Sigma_2$.  The strip map associated
to $(\Sigma_2,V_2)$ maps the squares in $\Sigma_1$
to the squares in $\Sigma_2$, translating along the diagonal
lines. This is shown in Figure 8.3.

\begin{center}
\resizebox{!}{4in}{\includegraphics{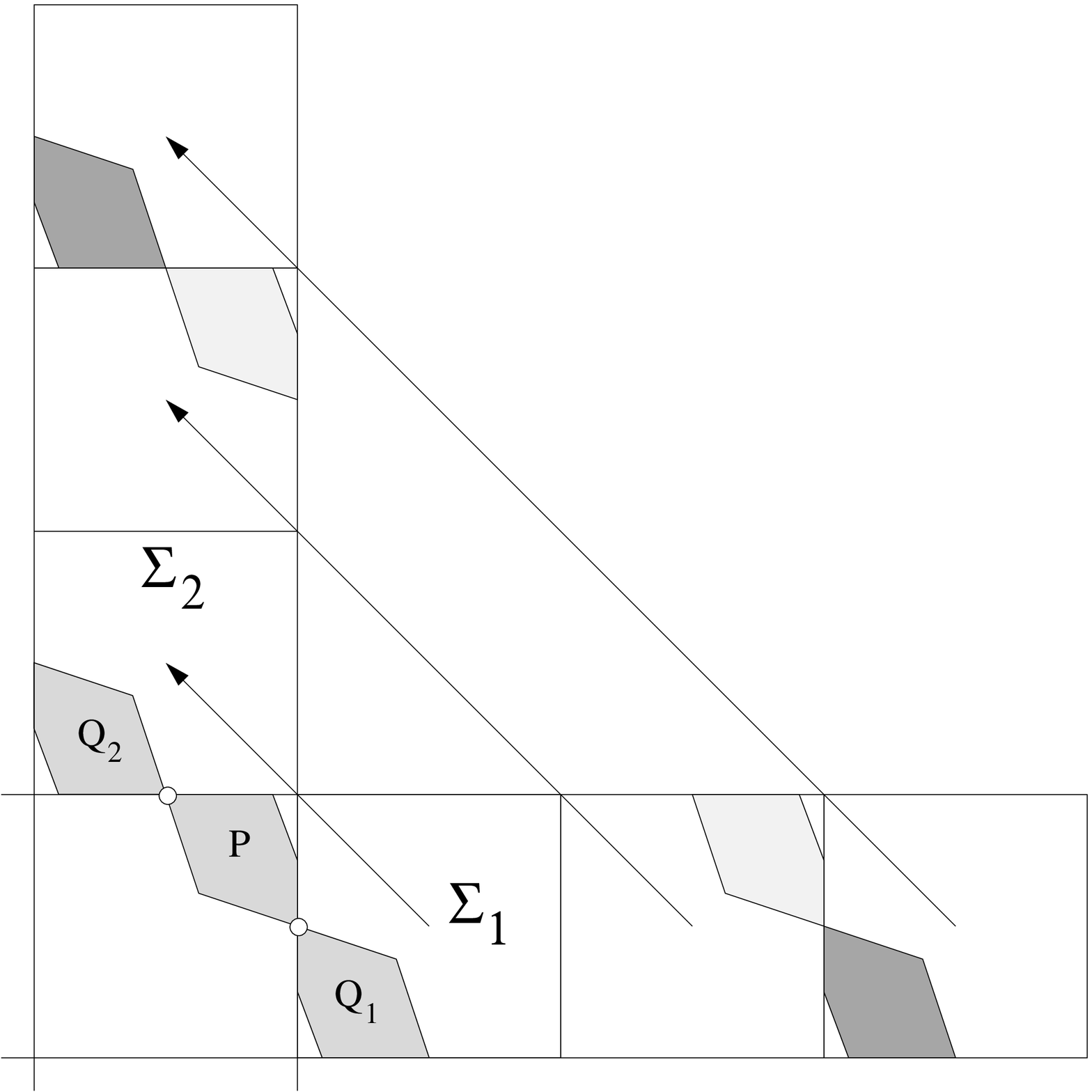}}
\newline
{\bf Figure 8.3:\/} Action of the strip map
\end{center}

Now we identify the unlabelled polygons in the picture.
\begin{itemize}

\item The light  polygon at bottom right is the projection to the plane of
$P_1^M$ for some integer $M$.  (The picture shows the case when $M=2$.)

\item The dark polygon at bottom right is the projection to the plane of
$Q_1^M$.

\item The light  polygon at top left is the projection to the plane of
$\psi(P_1^M)$.

\item The dark  polygon at top left is the projection to the plane of
$\psi(Q_1^M)$.

\end{itemize}

Note that $\psi(P_1)=P_1$ and $\psi(Q_1)=Q_2$.
From the picture, we can see that $\psi$ maps any point
in $\Sigma_1$ lying between $R_1$ and $R_1^M$ to 
a point in $\Sigma_2$ lying between $\psi(R_1)$ and
$\psi(R_1^M)$.      

Now, $\Sigma_2$ is also tiled by parallelograms which are
translates of the parallelogram $\Sigma_2 \cap \Sigma_3$.
Figure 8.4 shows the picture. 

\begin{center}
\resizebox{!}{4in}{\includegraphics{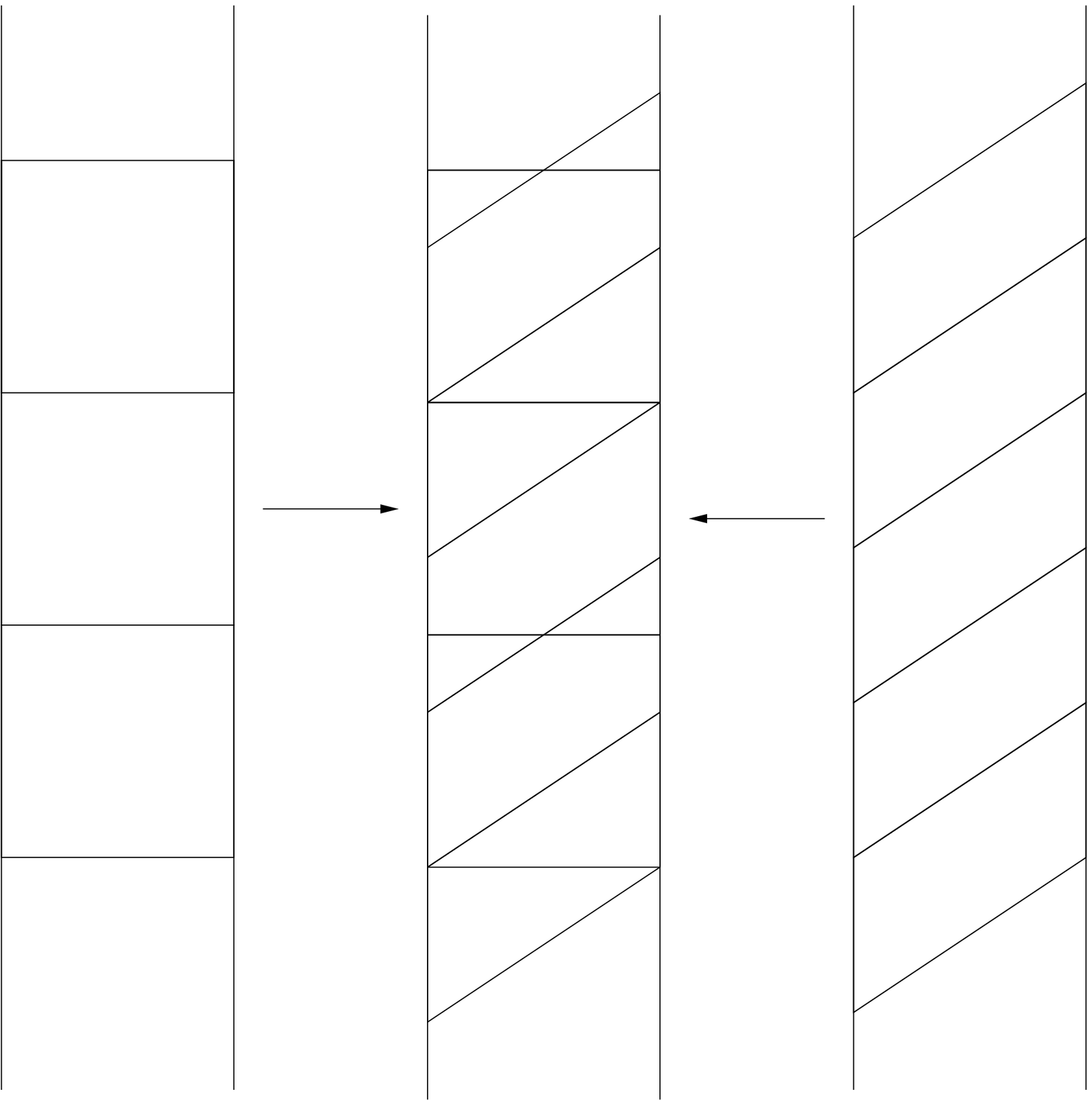}}
\newline
{\bf Figure 8.4:\/} Two superimposed tilings
\end{center}

Let $h$ denote the length of the vertical edge of
$\Sigma_2 \cap \Sigma_3$.   We have
\begin{equation}
h=\frac{A_3}{A_2}.
\end{equation}
From this we see that
\begin{equation}
\psi(P_1^M)=P_2^{M'}; \hskip 30 pt
\psi(Q_1^M)=Q_2^{M'}; \hskip 30 pt
M'=M \frac{A_2}{A_3}.
\end{equation}
provided that $M'$ is also an integer.

Now we use the fact that $P$ is quasirational.  We set
\begin{equation}
M=mD_1; \hskip 30 pt m \in \Z.
\end{equation}
This gives us
\begin{equation}
M'=mD_2 \in \Z.
\end{equation}

Observing that the same argument can be carried out for 
any choice of $j$ and not just $j=1$, we see that we have
completed the proof.
\endproof

Lemma \ref{quasi1} immediately implies that 
any point in $\Sigma_j$ between $R_j$ and
$R_j^{mD_j}$ has a forwards bounded orbit.
Hence, all orbits of $\psi$ are bounded.   By Corollary \ref{bijection2},
all forward outer billiards orbits are bounded.  By symmetry,
all backward outer billiards orbits are also bounded.
This completes the proof.    
\newline
\newline
{\bf Remarks:\/}
\newline
(i) The outer billiards orbits $O_m$ corresponding to
$R_j^{mD_j}$, at least for $m$ large,
are known as {\it necklace orbits\/} in [{\bf GS\/}].
These orbits consist of a ``necklace'' of polygons
winding once around $P$, and touching vertex to
vertex as shown in Figure 8.3.
\newline
(ii) We have stated the proof in such a way that it appears
to require the full force of the Pinwheel Theorem.  However,
this is not the case.  Outside of a large compact region,
there is an obvious correspondence between the pinwheel
maps and the outer billiards maps.  For $m$ large, our
argument shows that the region of the plane bounded
by the orbit $O_m$ is invariant under the outer billiards map.
This gives an exhaustion of $\R^2$ by bounded invariant sets.

\end{document}